\documentclass[a4j,11pt]{article}

\usepackage{amsmath}
\usepackage{amssymb}
\usepackage{amsthm}
\usepackage{array}
\usepackage{amscd}
\usepackage{cases}
\usepackage{bm}
\usepackage[all]{xy}
\usepackage{authblk}

\theoremstyle{plain}
\newtheorem{thm}{Theorem}[section]
\newtheorem{prop}[thm]{Proposition}

\theoremstyle{definition}
\newtheorem{df}{Definition}

\theoremstyle{definition}
\newtheorem{rem}{Remark}
\newtheorem{eg}[thm]{Example}

\theoremstyle{remark}

\numberwithin{equation}{section}

\title{Four-dimensional Painlev\'e-type equations associated with ramified linear equations I: Matrix Painlev\'e systems}
\author{Hiroshi KAWAKAMI\thanks{\texttt{kawakami@gem.aoyama.ac.jp}}}
\affil{College of Science and Engineering, Aoyama gakuin university,
5-10-1 Fuchinobe, Chuo-ku, Sagamihara-shi, Kanagawa 252-5258, Japan.}
\date{}

\begin{document}
\maketitle

\begin{abstract}
As a sequel to~\cite{KNS},
this series of three papers constructs the complete degeneration scheme of four-dimensional Painlev\'e-type equations
which includes the Painlev\'e-type equations associated with linear systems of ramified type.
In the present paper, we consider the degeneration of Painlev\'e-type equations which we call the matrix Painlev\'e systems.
\end{abstract}


\tableofcontents

\section{Introduction}\label{sec:intro}

The Painlev\'e equations, which were discovered by Painlev\'e~\cite{P} and Gambier~\cite{Gm} in the early twentieth century,
are non-linear second order ordinary differential equations that define new special functions.
They are closely related to the classical special functions (such as the Gauss hypergeometric function, the Bessel function, and so on)
and elliptic functions.

Originally, the Painlev\'e equations were classified into six equations.
We often denote them as $P_{\mathrm{I}},\ldots,P_{\mathrm{VI}}$.
However, from the viewpoint of the theory of initial value spaces (\cite{O1}),
it is natural to divide the third Painlev\'e equation into three types according to the values of its parameters.
We denote them as $P_{\mathrm{III}(D_6)}, P_{\mathrm{III}(D_7)}$, and $P_{\mathrm{III}(D_8)}$.
The so-called third Painlev\'e equation is then $P_{\mathrm{III}(D_6)}$.
We thus consider that there are eight types of Painlev\'e equations~(\cite{Sak1}).

As is well known, the Painlev\'e equations can be written in Hamiltonian form (\cite{Ok, OKSO}).
Here we give the eight Hamiltonians:
{\allowdisplaybreaks
\begin{align*} 
& t(t-1)H_{\rm VI}\left({\alpha , \beta\atop\gamma, \delta 
 };t;q,p\right)=\;q(q-1)(q-t)p^2\\ 
 & \hspace{13em}+\{ \delta q(q-1)-(2\alpha +\beta +\gamma +\delta )q(q-t)+\gamma
 (q-1)(q-t)\} p\\ 
 & \hspace{13em}+\alpha (\alpha +\beta )(q-t),\\ 
& tH_{\rm V}\left({\alpha , \beta \atop \gamma };t;q,p\right)=\;p(p+t)q(q-1) 
+\beta pq+\gamma p-(\alpha +\gamma )tq,\\ 
& H_{\rm IV}\left(\alpha , \beta;t;q,p\right)=\; 
pq(p-q-t)+\beta p+\alpha q,\quad
 tH_{\mathrm{III}(D_6)}\left(\alpha , \beta ;t;q,p\right)=\; 
p^2q^2-(q^2-\beta q-t)p-\alpha q,\\ 
& tH_{\mathrm{III}(D_7)}\left(\alpha;t;q,p\right)=\; 
p^2q^2+\alpha qp+tp+q,\quad
 tH_{\mathrm{III}(D_8)}\left(t;q,p\right)=\; 
p^2q^2+qp-q-\frac{t}{q},
\\ 
& H_{\rm II}\left(\alpha;t;q,p\right)=\; 
p^2-(q^2+t)p-\alpha q,\quad \hspace{3em}
 H_{\rm I}\left(t;q,p\right)=\; 
p^2-q^3-tq,
\end{align*}
}
where $H_{\mathrm{J}}$ is the Hamiltonian corresponding to $P_{\mathrm{J}}$.

An important aspect of the Painlev\'e equations for us is their relation to linear differential equations;
that is, \textit{isomonodromic deformations} of linear differential equations.
The first example was given by R. Fuchs~\cite{F} in the case of the sixth Painlev\'e equation.
An isomonodromic deformation is a deformation of a linear differential equation that does not change its ``monodromy data'', see \cite{JMU} for details.
We only give here a brief description of isomonodromic deformations.
Throughout this series of papers, linear differential equations are written in the form of first order systems.

Consider a system of first order linear differential equations:
\begin{equation}\label{eq:linear_system}
\frac{dY}{dx}=
A(x,u)Y
\end{equation}
where $A(x,u)$ is a matrix whose entries are rational functions in $x$ and depend on some parameters $u=(u_1,\ldots,u_n)$.


The isomonodromic deformation of (\ref{eq:linear_system}) is equivalent to an existence of matrices $B_i(x,u) \, (i=1,\ldots,n)$, rational in $x$,
such that
\begin{align*}
\left\{
\begin{aligned}
\frac{\partial Y}{\partial x}&=A(x,u)Y \\
\frac{\partial Y}{\partial u_i}&=B_i(x,u)Y \quad (i=1,\ldots,n)
\end{aligned}
\right.
\end{align*}
are completely integrable.
Then the isomonodromic deformation equations of (\ref{eq:linear_system}) are written as
\begin{align*}
&\frac{\partial A(x,u)}{\partial u_i}-\frac{\partial B_i(x,u)}{\partial x}+[A(x,u),B_i(x,u)]=O,\\
&\frac{\partial B_i(x,u)}{\partial u_j}-\frac{\partial B_j(x,u)}{\partial u_i}+[B_i(x,u),B_j(x,u)]=O.
\end{align*}
Hereafter we use the term {\it Painlev\'e-type equations} instead of isomonodromic deformation equations.

For example,
suppose that the system (\ref{eq:linear_system}) has the following form
\begin{equation}\label{eq:Fuchs}
\frac{dY}{dx}=
\sum_{i=1}^n \frac{A_i}{x-u_i}Y \quad (A_i \in M_m(\mathbb{C}))
\end{equation}
with some generic conditions.
Then $B_i$'s are known to be given by $B_i(x,u)=-\frac{A_i}{x-u_i}$~(\cite{Sc}).
%



Here we review the classification of the Painlev\'e equations in terms of isomonodromic deformations.
Among the Painlev\'e equations, the sixth Painlv\'e equation is the ``source equation'' in the sense that
all the other Painlev\'e equations can be obtained from $P_{\mathrm{VI}}$ through degenerations.
First, we look at the degeneration scheme of Painlev\'e equations, based on the original classification
which classifies the Painlev\'e equations into six types.
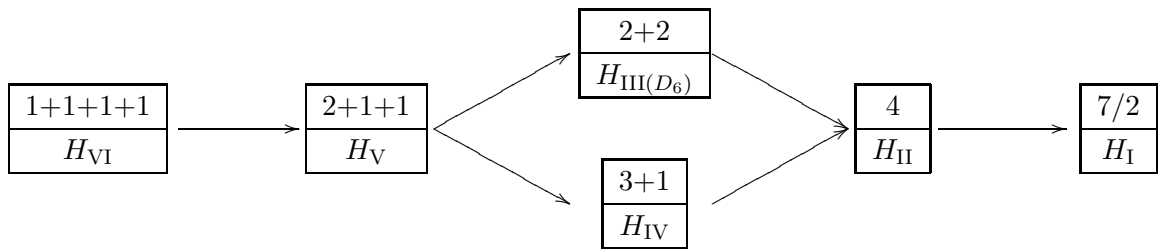
\vspace{3mm}


\begin{xy}
{(3,0) *{\begin{tabular}{|c|}
\hline
1+1+1+1\\
\hline
$H_{\rm VI}$\\
\hline
\end{tabular}
}},
{(40,0) *{\begin{tabular}{|c|}
\hline
2+1+1\\
\hline
$H_{\rm V}$\\
\hline
\end{tabular}
}},
{\ar (15,0);(31,0)},
{\ar (49,0);(67,10)},
{\ar (49,0);(67,-10)},
{(77,10) *{\begin{tabular}{|c|}
\hline
2+2\\
\hline
$H_{\mathrm{III}(D_6)}$\\
\hline
\end{tabular}}},
{\ar (86,10);(104,0)},
{\ar (86,-10);(104,0)},
{(77,-10) *{\begin{tabular}{|c|}
\hline
3+1\\
\hline
$H_{\rm IV}$\\
\hline
\end{tabular}}},
{(110,0) *{\begin{tabular}{|c|}
\hline
4\\
\hline
$H_{\rm II}$\\
\hline
\end{tabular}}},
{\ar (116,0);(133,0)},
{(140,0) *{\begin{tabular}{|c|}
\hline
7/2\\
\hline
$H_{\mathrm{I}}$\\
\hline
\end{tabular}}},
\end{xy}

\vspace{3mm}

The number in the upper half of each box is the ``singularity pattern'', 
which expresses Poincar\'e ranks of the singularities of the associated linear system.
More specifically, each number partitioned by + expresses the number ``Poincar\'e rank +1'' (see Section~\ref{sec:HTL}).
In particular, 1+1+1+1 means that corresponding linear system has four regular singular points (thus it is a Fuchsian system).

However, in the above scheme, the third Painlev\'e equations of type $D_7^{(1)}, D_8^{(1)}$ are absent.
Moreover, in view of associated linear systems, 
the degeneration $H_{\mathrm{II}} \to H_{\mathrm{I}}$ has different meaning from the other degenerations.
Namely, the other degenerations correspond to confluences of singular points of associated linear systems,
while the degeneration $H_{\mathrm{II}} \to H_{\mathrm{I}}$ corresponds to the degeneration of the Jordan canonical form
of the leading coefficient of (or, more precisely, the degeneration of the HTL canonical form at) the irregular singular point.

Considering the other possible degenerations of HTL canonical forms and confluences of singular points,
one can obtain the following ``complete'' degeneration scheme~\cite{OO}:
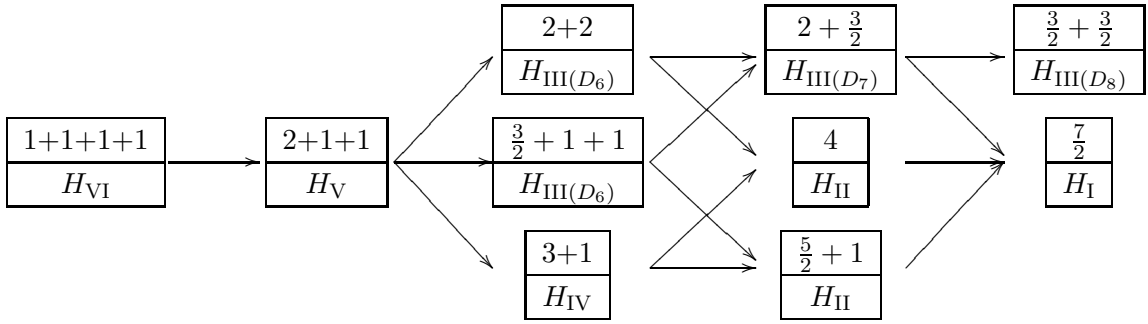
\vspace{3mm}

\begin{xy}
{(3,0) *{\begin{tabular}{|c|}
\hline
1+1+1+1\\
\hline
$H_{\mathrm{VI}}$\\
\hline
\end{tabular}
}},
{(35,0) *{\begin{tabular}{|c|}
\hline
2+1+1\\
\hline
$H_{\mathrm{V}}$\\
\hline
\end{tabular}
}},
{\ar (14,0);(26,0)},
{\ar (44,0);(57,14)},
{\ar (44,0);(57,0)},
{\ar (44,0);(57,-14)},
{(67,15) *{\begin{tabular}{|c|}
\hline
2+2\\
\hline
$H_{\mathrm{III}(D_6)}$\\
\hline
\end{tabular}}},
{(67,0) *{\begin{tabular}{|c|}
\hline
$\frac{3}{2}+1+1$\\
\hline
$H_{\mathrm{III}(D_6)}$\\
\hline
\end{tabular}}},
{\ar (78,14);(92,14)},
{\ar (78,14);(92,1)},
{\ar (78,0);(92,13)},
{\ar (78,0);(92,-13)},
{\ar (78,-14);(92,-1)},
{\ar (78,-14);(92,-14)},
{(67,-15) *{\begin{tabular}{|c|}
\hline
3+1\\
\hline
$H_{\mathrm{IV}}$\\
\hline
\end{tabular}}},
{(102,15) *{\begin{tabular}{|c|}
\hline
$2+\frac{3}{2}$\\
\hline
$H_{\mathrm{III}(D_7)}$\\
\hline
\end{tabular}}},
{(102,0) *{\begin{tabular}{|c|}
\hline
4\\
\hline
$H_{\mathrm{II}}$\\
\hline
\end{tabular}}},
{\ar (112,14);(125,1)},
{\ar (112,14);(125,14)},
{\ar (112,0);(125,0)},
{\ar (112,-14);(125,0)},
{(102,-15) *{\begin{tabular}{|c|}
\hline
$\frac{5}{2}+1$\\
\hline
$H_{\mathrm{II}}$\\
\hline
\end{tabular}}},
{(135,15) *{\begin{tabular}{|c|}
\hline
$\frac{3}{2}+\frac{3}{2}$\\
\hline
$H_{\mathrm{III}(D_8)}$\\
\hline
\end{tabular}}},
{(135,0) *{\begin{tabular}{|c|}
\hline
$\frac{7}{2}$\\
\hline
$H_{\mathrm{I}}$\\
\hline
\end{tabular}}},
\end{xy}

\vspace{3mm}
\begin{rem}
In the case of  the Painlev\'e equations,
since (standard) corresponding linear systems are rank 2,
the degenerations of HTL canonical forms are caused by the degenerations of Jordan canonical forms.
However, in general, a degeneration of an HTL canonical form does not always correspond to a degeneration of a Jordan canonical form
(this will be shown in the forthcoming papers~\cite{K2, K3}).
\qed
\end{rem}

Recently, there have been many generalizations of the Painlev\'e equations in the literature.
What is important here is that they can be written as compatibility conditions of linear differential equations.
Thus we can say that they describe isomonodromic deformations of some linear differential equations.
The purpose of the present series of papers is to classify those of four-dimensional phase space using associated linear systems.
More specifically, 
we construct the degeneration scheme starting from suitable Fuchsian systems
so that we obtain the classification of four-dimensional Painlev\'e-type equations.

What is the ``source equation'' of the degeneration scheme of four-dimensional Painlev\'e-type equations?
The recent development of the theory of Fuchsian systems (\cite{Katz, Os, HF}) enables us to
answer the question.
That is, there are four Fuchsian systems (thus corresponding four Painlev\'e-type equations~\cite{Sak2}) which should be placed 
at the starting points of the degeneration scheme of four-dimensional Painlev\'e type equations:
one admits two-dimenional deformation and corresponds to the Garnier system~\cite{G},
the remaining three correspond to the so-called Fuji-Suzuki system~\cite{FS1, Ts} (abbreviated to FS-system),
the Sasano system~\cite{Ss}, and the sixth matrix Painlev\'e system~\cite{B2, K}.

In \cite{KNS}, the authors constructed the degeneration scheme starting with the above four Fuchsian systems
by considering confluences of singular points.
As the result, they obtained the degeneration scheme of four-dimensional Painlev\'e-type equations
associated with unramified linear equations.
Note that the degeneration scheme of the Garnier system in two variables had already been obtained by Kimura~\cite{Ki}.

However, the confluence of singularities of linear equations is not sufficient to produce all the four-dimensional Painlev\'e-type equations.
The aim of this series of papers is, by considering the degeneration of HTL canonical forms,
to obtain the ``complete'' degeneration scheme of four-dimensional Painlev\'e-type equations.


In the present paper, we focus on degenerations of HTL canonical forms starting from the Fuchsian system 
corresponding to what we call the sixth matrix Painlev\'e system.

This paper is organized as follows.
In Section~\ref{sec:HTL}, we explain HTL canonical forms at singular points of linear systems
and demonstrate the degeneration of an HTL form.
In Section~\ref{sec:Lax_pairs}, we present the Hamiltonians 
associated with ramified linear systems together with their Lax pairs.
Section~\ref{sec:Laplace} is devoted to describe correspondences through the Laplace transform.

We give in advance the degeneration scheme of the matrix Painlev\'e systems.

\vspace{2mm}

{\scriptsize
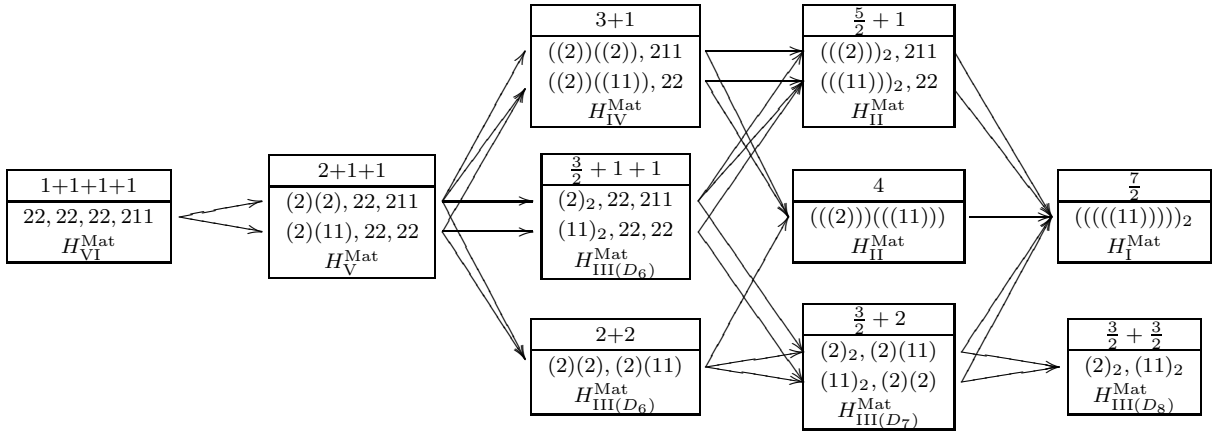
\begin{xy}
{(0,0) *{\begin{tabular}{|c|}
\hline
1+1+1+1\\
\hline
$22,22,22,211$\\
$H^{\mathrm{Mat}}_{\mathrm{VI}}$\\
\hline
\end{tabular}
}},
{\ar (12,0);(23,2)},
{\ar (12,0);(23,-2)},
{(35,0) *{\begin{tabular}{|c|}
\hline
2+1+1\\
\hline
$(2)(2),22,211$\\
$(2)(11),22,22$\\
$H^{\mathrm{Mat}}_{\mathrm{V}}$\\
\hline
\end{tabular}
}},
{\ar (47,2);(58,22)},
{\ar (47,2);(58,17)},
{\ar (47,2);(58,-19)},
{\ar (47,2);(59,2)},
{\ar (47,-2);(58,17)},
{\ar (47,-2);(59,-2)},
{\ar (47,-2);(58,-19)},
{(70,20) *{\begin{tabular}{|c|}
\hline
3+1\\
\hline
$((2))((2)),211$\\
$((2))((11)),22$\\
$H^{\mathrm{Mat}}_{\mathrm{IV}}$\\
\hline
\end{tabular}}},
{\ar (82,22);(93,0)},
{\ar (82,22);(95,22)},
{\ar (82,18);(93,0)},
{\ar (82,18);(95,18)},
{(70,0) *{\begin{tabular}{|c|}
\hline
$\frac{3}{2}+1+1$\\
\hline
$(2)_2,22,211$\\
$(11)_2,22,22$\\
$H^{\mathrm{Mat}}_{\mathrm{III}(D_6)}$\\
\hline
\end{tabular}}},
{\ar (81,2);(95,22)},
{\ar (81,2);(95,18)},
{\ar (81,2);(95,-18)},
{\ar (81,-2);(95,18)},
{\ar (81,-2);(95,-22)},
{(70,-20) *{\begin{tabular}{|c|}
\hline
2+2\\
\hline
$(2)(2),(2)(11)$\\
$H^{\mathrm{Mat}}_{\mathrm{III}(D_6)}$\\
\hline
\end{tabular}}},
{\ar (82,-20);(93,0)},
{\ar (82,-20);(95,-18)},
{\ar (82,-20);(95,-22)},
{(105,20) *{\begin{tabular}{|c|}
\hline
$\frac52+1$ \\
\hline
$(((2)))_2, 211$\\
$(((11)))_2, 22$\\
$H^{\mathrm{Mat}}_{\mathrm{II}}$\\
\hline
\end{tabular}}},
{\ar (115,22);(128,0)},
{\ar (115,17);(128,0)},
{(105,0) *{\begin{tabular}{|c|}
\hline
4\\
\hline
$(((2)))(((11)))$\\
$H^{\mathrm{Mat}}_{\mathrm{II}}$\\
\hline
\end{tabular}}},
{\ar (117,0);(128,0)},
{(105,-20) *{\begin{tabular}{|c|}
\hline
$\frac32+2$\\
\hline
$(2)_2, (2)(11)$\\
$(11)_2, (2)(2)$\\
$H^{\mathrm{Mat}}_{\mathrm{III}(D_7)}$\\
\hline
\end{tabular}}},
{\ar (116,-18);(129,-20)},
{\ar (116,-18);(128,0)},
{\ar (116,-22);(129,-20)},
{\ar (116,-22);(128,0)},
{(139,0) *{\begin{tabular}{|c|}
\hline
$\frac72$ \\
\hline
$(((((11)))))_2$\\
$H^{\mathrm{Mat}}_{\mathrm{I}}$\\
\hline
\end{tabular}}},
{(139,-20) *{\begin{tabular}{|c|}
\hline
$\frac32+\frac32$ \\
\hline
$(2)_2, (11)_2$\\
$H^{\mathrm{Mat}}_{\mathrm{III}(D_8)}$\\
\hline
\end{tabular}}}
\end{xy}
}

\vspace{2mm}

Degeneration schemes of the Sasano systems, the FS systems, and the Garnier systems will appear in
our forthcoming papers~\cite{K2,K3}.

\bigskip 

\noindent 
\textbf{Acknowledgements} 

\noindent 
The author wishes to thank Professors Hidetaka Sakai and Akane Nakamura
for their invaluable advices and comments.
The author also thanks Professor Kazuki Hiroe.
The notation of spectral types for ramified singularities in this series of papers follows his suggestion.

\section{HTL canonical forms and their degeneration}\label{sec:HTL}
As mentioned earlier, to classify the four-dimensional Painlev\'e-type equations,
we attach a linear system to each Painlev\'e-type equation.
Hence we would like to introduce simple notations of linear systems; 
that is, Riemann schemes and, their abbreviated symbols, spectral types.

The Riemann scheme and the spectral type of a linear system will be defined through the collection of ``canonical forms'' of
the linear system.
Here the canonical form is determined at each singular point.
We call the canonical form the HTL canonical form (or HTL form) of the system.

\subsection{HTL canonical forms}\label{sec:HTL}
Here we recall the HTL canonical form of a linear system at a singular point.
Note that for a linear system
\[
\frac{dY}{dx}=A(x)Y,
\]
a transformation of the dependent variable $Y=P(x)Z$ by an invertible matrix $P(x)$ yields the
following system:
\[
\frac{dZ}{dx}=\left( P(x)^{-1}A(x)P(x)-P(x)^{-1}\frac{dP(x)}{dx} \right)Z.
\]
We write the coefficient matrix $P^{-1}AP-P^{-1}\frac{dP}{dx}$ as $A^P(x)$.
This kind of transformation $A(x) \mapsto A^P(x)$ is called a gauge transformation.
We often express the above transformation of the dependent variable as $Y \to PY$
(i.e. we write the new dependent variable as $Y$ again).

Linear systems that we treat in this series of papers are those with rational function coefficients.
Such a system is generally written as follows:
\begin{equation}
\frac{dY}{dx}=
\left(\sum_{\nu=1}^n\sum_{k=0}^{r_{\nu}}\frac{A_{\nu}^{(k)}}{(x-u_{\nu})^{k+1}}
+\sum_{k=1}^{r_{\infty}}A_{\infty}^{(k)}x^{k-1}
\right)Y, \quad
A_j^{(k)} \in M_m(\mathbb{C}).
\end{equation}
This system has singular points at $x=u_1,\ldots,u_n$, and $\infty$.
Choosing one of the singular points and taking the local coordinates $z=x-u_\nu$ or $z=1/x$
according to the choice,
we write the system at $z=0$ as follows:
\begin{equation}\label{eq:Laurent}
\frac{dY}{dz}=\left(
\frac{A_0}{z^{r+1}}+\frac{A_1}{z^{r}}+\cdots+A_{r+1}+A_{r+2} z+\cdots
\right)Y.
\end{equation}
We denote the field of formal Laurent series in $z$ by $\mathbb{C}(\!(z)\!)$,
and the field of Puiseux series $\cup_{p > 0}\mathbb{C}(\!(z^{\frac1p})\!)$ by $\mathcal{K}_z$.

We can convert (\ref{eq:Laurent}) to the HTL canonical form.
Here the definition of the HTL canonical form is as follows.
\begin{df}
An element in $M_m(\mathcal{K}_z)$ of the form
\begin{equation}
\frac{D_0}{z^{l_0}}+\frac{D_1}{z^{l_1}}+\cdots+\frac{D_{s-1}}{z^{l_{s-1}}}+\frac{\Theta}{z^{l_s}}
\end{equation}
where
	\begin{itemize}
		\item $l_0, \ldots, l_{s}\ (l_0 > l_1 > \cdots > l_{s-1}>l_s = 1)$ are rational numbers,
		
		\item $D_0, \ldots, D_{s-1}$ are diagonal matrices,
		
		\item $\Theta$ is a (not necessarily diagonal) Jordan matrix which commutes with all $D_j$'s,
	\end{itemize}
is called an {\it HTL canonical form} or {\it HTL form} for short.
\qed
\end{df}
The following theorem is fundamental for us.
\begin{thm}[Hukuhara \cite{Huk}, Turrittin \cite{Tur}, Levelt \cite{Lev}]
For any
\begin{equation}\label{eq:Laurent2}
A(z)=\frac{A_0}{z^{r+1}}+\frac{A_1}{z^{r}}+\cdots+A_{r+1}+A_{r+2} z+\cdots, \quad
A_j \in M_m(\mathbb{C}),
\end{equation}
there exists a matrix $P(z) \in \mathrm{GL}_m(\mathcal{K}_z)$ such that $A^P(z)$ is an HTL form
\begin{equation}\label{eq:HTLform}
A^P(z)=
\frac{D_0}{z^{l_0}}+\frac{D_1}{z^{l_1}}+\cdots+\frac{D_{s-1}}{z^{l_{s-1}}}+\frac{\Theta}{z^{l_s}}.
\end{equation}
%
%
Here $l_0, \ldots, l_{s}$ are uniquely determined only by the original system (\ref{eq:Laurent2}).

If the following
\begin{equation}
\frac{\tilde{D}_0}{z^{l_0}}+\frac{\tilde{D}_1}{z^{l_1}}+\cdots+\frac{\tilde{D}_{s-1}}{z^{l_{s-1}}}
+\frac{\tilde{\Theta}}{z^{l_s}}
\end{equation}
is another HTL canonical form corresponding to the same system~(\ref{eq:Laurent2}),
there exists a constant matrix $g \in \mathrm{GL}_m(\mathbb{C})$ and a natural number $k \in \mathbb{Z}_{\ge 1}$ such that
\begin{equation}
\tilde{D}_j=g^{-1}D_j g, \quad \exp(2\pi i k \tilde{\Theta})=g^{-1}\exp(2\pi i k \Theta)g
\end{equation}
hold.
\end{thm}
We call (\ref{eq:HTLform}) the HTL canonical form (or HTL form) of (\ref{eq:Laurent2}).
The number $l_0-1$ is called the {\it Poincar\'e rank} of the singular point.
When there is a rational number $l_j$ that is not an integer, the singular point is called a {\it ramified} irregular singular point.
A linear system is said to be of \textit{ramified type} if it has a ramified irregular singular point.

When we express the Poincar\'e ranks of the singular points of a given system,
we attach the number ``Poincar\'e rank +1'' to each singular point and
connect them with ``+''.
We call it the {\it singularity pattern} of the system.

In this series of papers, the residue matrix $\Theta$ of an HTL form is always diagonal (see Section~\ref{sec:shearing}).
In the next subsection, we define Riemann schemes and spectral types for such cases.

\subsection{Riemann schemes and spectral types}
\subsubsection{Unramified singularities and refining sequences of partitions}\label{sec:RSP}
The Riemann scheme of a linear system is defined to be the table of HTL forms of the system at all singular points.
Here we introduce a special notation for HTL forms.

For an $m \times m$ HTL form
\begin{equation}\label{eq:HTL_form}
\frac{D_0}{z^{l_0}}+\frac{D_1}{z^{l_1}}+\cdots+\frac{D_{s-1}}{z^{l_{s-1}}}+\frac{\Theta}{z},
\end{equation}
let $d \in \mathbb{Z}_{>0}$ be the minimum element of 
$\{ k \in \mathbb{Z}_{>0} \,|\,  k l_j \in \mathbb{Z} \ (j=0,\ldots,s-1) \}$.
Then the canonical form (\ref{eq:HTL_form}) is rewritten in the following form
\begin{equation}
\frac{T_{0}}{z^{\frac{b}{d}+1}}+\frac{T_{1}}{z^{\frac{b-1}{d}+1}}+\cdots+\frac{T_{b-1}}{z^{\frac{1}{d}+1}}+\frac{\Theta}{z}
\end{equation}
where some of $T_j$'s may be the zero matrix.
Then, by writing the diagonal entries of $T_j$'s and $\Theta$ as $t^j_k$ and $\theta_k\ (k=1, \ldots, m)$
respectively, we can express the canonical form as follows:
\begin{equation}
\begin{array}{c}
x=u_i \ \left( \frac{1}{d} \right) \\
\overbrace{\begin{array}{ccccc}
	t^0_1  & t^1_1 & \ldots & t^{b-1}_1  & \theta_1\\
	\vdots & \vdots    &        & \vdots & \vdots \\
	t^0_m  & t^1_m & \ldots & t^{b-1}_m  & \theta_m
	\end{array}}
\end{array}.
\end{equation}
When the position of a singular point under consideration is $\infty$, then of course we write ``$x=\infty$''.
The symbol $(\frac1d)$ is omitted when $d=1$. 
The table of the above expressions of HTL forms for all singularities of a linear system is called the {\it Riemann scheme} of the system.

It is shown in \cite{KNS} that the feature of an HTL form of a linear system at an \textit{unramified} singular point
is well described by the refining sequence of partitions.
\begin{df}\label{def:rsp}
Let $\lambda=\lambda_1\ldots \lambda_p$ and $\mu=\mu_1\ldots \mu_q$ be partitions of a natural number $m$:
$\lambda_1+\cdots+\lambda_p=\mu_1+\cdots+\mu_q=m$．
Here we assume that $\lambda_i$'s and $\mu_i$'s are not necessarily arranged in descending or ascending order.

If there exist a disjoint decomposition $\{ 1,2,\ldots , p\}=I_1\coprod \cdots \coprod I_q$ of the index set of $\lambda$ such that $\mu_k=\sum_{j\in I_k}\lambda_j$ holds,
then we call $\lambda$ a {\it refinement} of $\mu$.

Let $[p_0,\ldots,p_r]$ be an $(r+1)$-tuple of partitions of $m$.
When $p_{i+1}$ is a refinement of $p_i$ for all $i\ (i=0,\ldots,r-1)$, we call $[p_0,\ldots,p_r]$ a {\it refining sequence of partitions} or {\it RSP} for short.
\qed
\end{df}
%
%
The appearance of the RSP is a  consequence of the following proposition.

\begin{prop}{(block diagonalization)}\label{thm:block_diag}
For any formal Laurent series
\begin{equation}\label{eq:irreg_sys}
A(z)=\frac{1}{z^{r+1}}\left( A_0+A_1z+\cdots \right)
\end{equation}
where the eigenvalues of $A_0$ are $\lambda_1,\ldots,\lambda_n$ and their multiplicities $m_1,\ldots,m_n$ respectively,
we can choose a formal power series $P(z)$ so that the gauge transformation by $P(z)$ changes {\rm (\ref{eq:irreg_sys})} to
the following form:
\begin{align}
&A^P(z)
=
\begin{pmatrix}
B_1 & & \\
& \ddots & \\
& & B_n
\end{pmatrix}, \quad
B_k=\frac{1}{z^{r+1}}\left( B^k_0+B^k_1z+\cdots \right)
\end{align}
where $B^k_0=\lambda_k I_{m_k}+N_k$ and $N_k$ is nilpotent.
\end{prop}
Therefore the system corresponding to (\ref{eq:irreg_sys}) is formally reduced to the following $n$ systems:
\begin{equation}\label{eq:dec_system}
\frac{dZ_k}{dz}=B_k Z_k \quad (k=1,\ldots,n).
\end{equation}

To compute the HTL form for a given linear system, we first decompose the linear system according to
Proposition~\ref{thm:block_diag}.
If the leading matrices $B^1_0, \ldots, B^n_0$ are semisimple (i.e. if they are scalar matrices),
we focus on the next matrices $B^1_1, \ldots, B^n_1$.
If they are all semisimple, then (\ref{eq:dec_system}) again decompose along each eigenvalue of $B^k_1\ (k=1,\ldots,n)$.
In this way, the RSP structure arises.
\begin{eg}
Consider the following HTL form:
{\small
\begin{align*}
\begin{pmatrix}
t^0_1 & & & \\
 & t^0_1 & & \\
 & & t^0_2 & \\
 & & & t^0_2
\end{pmatrix}
\frac{1}{z^3}+
\begin{pmatrix}
t^1_1 & & & \\
 & t^1_1 & & \\
 & & t^1_2 & \\
 & & & t^1_3
\end{pmatrix}
\frac{1}{z^2}+
\begin{pmatrix}
\theta_1 & & & \\
 & \theta_2 & & \\
 & & \theta_3 & \\
 & & & \theta_4
\end{pmatrix}
\frac{1}{z}.
\end{align*}
}
To this HTL form we attach the RSP $[22,211,1111]$, which is the spectral type of this HTL form.
The spectral type is also represented as $((11))((1)(1))$. We use the latter notation
%
%
(for detail, see~\cite{KNS}).
\qed
\end{eg}

\subsubsection{Spectral types of ramified singularities}\label{sec:shearing}
As we have seen, the spectral type of an \textit{unramified} singularity is an RSP.
If there appear non-semisimple matrices in the course of the above successive block diagonalizations, 
we need to perform the so-called shearing transformation.
Then the singularity is in general ramified.
As to a ramified singularity, its spectral type consists of  ``copies" of RSPs.

Here we briefly explain the shearing transformation.
According to Proposition~\ref{thm:block_diag}, it suffices to consider the reduction of the following system
\begin{equation}\label{eq:nilpotent_leading}
\frac{dY}{dz}=\frac{1}{z^{r+1}}\left( A_0+A_1z+\cdots+A_r z^r+\cdots \right)Y
\end{equation}
where $A_0$ is a Jordan matrix with only one eigenvalue to its HTL form.
By means of a scalar gauge transformation, we can shift $A_0$ by a scalar matrix.
Thus, without loss of generality, we can assume that $A_0$ is nilpotent.

The shearing transformation is a gauge transformation by a diagonal matrix, which is typically of the form
\begin{equation}\label{eq:shear_mat_general}
\mathrm{diag}(1, z^s, \ldots, z^{(m-1)s})
\end{equation}
where $s$ is a positive rational number.
By applying (more than once in general) the transformation of this kind with suitable value of $s$ to (\ref{eq:nilpotent_leading}),
we can make the leading matrix diagonalizable.
Here we illustrate with an example how the shearing transformation works.
For the general case, see for example \cite{BV, Huk, Wa}.
\begin{eg}
Consider the following system:
\begin{equation}\label{eq:eg_shear}
\frac{dY}{dx}=\left( \frac{A_0^{(-1)}}{x^2}+\frac{A_0^{(0)}}{x}+A_{\infty} \right)Y.
\end{equation}
Here
$A_0^{(-1)}$, $A_0^{(0)}$, and $A_{\infty}$ are given as follows:
\begin{align*}
A_0^{(-1)}&=t
\begin{pmatrix}
O_2 \\
I_2
\end{pmatrix}
\begin{pmatrix}
P &  I_2
\end{pmatrix},\quad
A_0^{(0)}=
\begin{pmatrix}
QP & Q \\
I_2 & -PQ+\theta^0
\end{pmatrix},\quad
A_{\infty}=
\begin{pmatrix}
O & I_2 \\
O & O
\end{pmatrix}
\end{align*}
where the matrices $Q$, $P$, and $\Theta$ are given by
\begin{align*}
Q=
\begin{pmatrix}
q_1 & u\\
-q_2/u & q_1
\end{pmatrix},
\quad
P=
\begin{pmatrix}
p_1/2 & -p_2 u\\
(p_2q_2-\theta^0-\theta^\infty_2)/u & p_1/2
\end{pmatrix},\quad
\Theta=
\begin{pmatrix}
\theta^{\infty}_2 & \\
 & \theta^{\infty}_3
\end{pmatrix}.
\end{align*}
Hereafter, we often abbreviate a scalar matrix $k I$ to $k$.

This system has an irregular singular point at $x=\infty$.
By changing the independent variable as $z=1/x$,
we rewrite the system (\ref{eq:eg_shear}) as follows:
\begin{align*}
\frac{dY}{dz}&=-
\left(
\frac{A_\infty}{z^2}+\frac{A_0^{(0)}}{z}+A_0^{(-1)}
\right)Y\\
&=:A(z)Y.
\end{align*}

In this case, a shearing matrix can be chosen as follows (though it has different form from (\ref{eq:shear_mat_general}))
\begin{equation}\label{eq:shear_mat}
S=
\begin{pmatrix}
I_2 & O \\
O & z^{1/2} I_2
\end{pmatrix}.
\end{equation}
Then we have
\begin{align*}
A^S(z)=
\begin{pmatrix}
O & -I \\
-I & O
\end{pmatrix}
\frac{1}{z^{3/2}}+
\begin{pmatrix}
-QP & O \\
O & PQ-\theta^0-1/2
\end{pmatrix}
\frac{1}{z}+
\begin{pmatrix}
O & -Q \\
-tP & O
\end{pmatrix}
\frac{1}{z^{1/2}}+
\begin{pmatrix}
O & O \\
O & -tI
\end{pmatrix}.
\end{align*}
Note that the leading term has semisimple coefficients.
In fact, let $G$ be
\begin{equation*}
G=
\begin{pmatrix}
-I & I \\
I & I
\end{pmatrix},
\end{equation*}
then we have
\begin{equation*}
A^{SG}(z)=
\begin{pmatrix}
I & O \\
O & -I
\end{pmatrix}
\frac{1}{z^{3/2}}+
\begin{pmatrix}
\Theta/2-1/4 & \frac{QP+PQ}{2}-\frac{\theta^0}{2}-\frac14 \\
\frac{QP+PQ}{2}-\frac{\theta^0}{2}-\frac14 & \Theta/2-1/4
\end{pmatrix}
\frac{1}{z}+\cdots.
\end{equation*}
Further, applying the gauge transformation by $z^{-1/4}$, we can eliminate $\frac{-\frac{1}{4}I_4}{z}$.
Hence the HTL form is
\begin{equation*}
\begin{pmatrix}
I & O \\
O & -I
\end{pmatrix}
\frac{1}{z^{3/2}}+
\begin{pmatrix}
\Theta/2 & O \\
O & \Theta/2
\end{pmatrix}
\frac{1}{z}.
\end{equation*}
This can be seen as the direct sum of $\frac{I_2}{z^{3/2}}+\frac{\Theta/2}{z}$ (whose spectral type is (11)) and its ``copy'' made by the action $z^{1/2} \mapsto -z^{1/2}$.
We write this spectral type as $(11)_2$.
Together with the HTL form at $x=0$ (it is not difficult to see),
we write the spectral type of the system (\ref{eq:eg_shear}) as $(11)_2,(2)(2)$.
The Riemann scheme is given in Section~\ref{sec:(11)_2,(2)(2)}.
\qed
\end{eg}
\begin{rem}
The matrix~(\ref{eq:shear_mat}) is applicable to the other ramified systems in this paper.
\qed
\end{rem}

Here 
we describe the structure of the HTL form (at a ramified irregular singularity) in a more general setting.
Let
\begin{equation}\label{eq:direct_summand}
T=\frac{T_{0}}{z^{\frac{b}{d}+1}}+\frac{T_{1}}{z^{\frac{b-1}{d}+1}}+\cdots+\frac{T_{b-1}}{z^{\frac{1}{d}+1}}+\frac{\Theta}{z}
\end{equation}
be the $d$-reduced (i.e. $0 \le \Re(\theta) < 1/d$ holds for any eigenvalue $\theta$ of $\Theta$) HTL form of $A(z) \in M_m(\mathbb{C}(\!(z)\!))$.
We set
\[
T_{\mathrm{irr}}=T-\frac{\Theta}{z}.
\]

Let $\Sigma \subset \mathbb{C}(\!(z^{1/d})\!)$ be the set of the eigenvalues of $T_{\mathrm{irr}} \in \mathrm{End}(V)$ where
$V=\mathbb{C}(\!(z^{1/d})\!) \otimes_{\mathbb{C}} \mathbb{C}^m$.
Clearly we have $V=\bigoplus_{p \in \Sigma} V(p)$
where $V(p)$ is the eigenspace of $T_{\mathrm{irr}}$ corresponding to $p \in \Sigma$.

Let $C_d$ be the (multiplicative) cyclic group: $C_d=\{ {\zeta_d}^j \,|\, j=0,\ldots, d-1\}$ where $\zeta_d=e^{2\pi i/d}$.
Note that $g \in C_d$ acts on $\mathbb{C}(\!(z^{1/d})\!)$ as
\[
f=\sum_j f_j (z^{1/d})^j \mapsto g \cdot f=\sum_j f_j g^{-j}(z^{1/d})^j.
\]
Then the following holds.
\begin{prop}[\cite{BV}]
The set $\Sigma$ is stable under the action of $C_d$,
and there exists a representation $\rho : C_d \to \mathrm{GL}(V)$
such that
\[
\rho(g)\Theta=\Theta \rho(g),\quad
\rho(g)(V(p))=V(g\cdot p) \quad
(g \in C_d, \, p \in \Sigma).
\]
\end{prop}

We note that there exists a direct sum decomposition of $\Sigma$:
\begin{equation*}
\Sigma=\Sigma_{d_1} \sqcup \cdots \sqcup \Sigma_{d_r}
\end{equation*}
where
\begin{equation*}
p \in \Sigma_{d_j} \stackrel{\mathrm{def}}{\iff}
\min\{ k \in \mathbb{Z}_{\ge 1} \, | \, {\zeta_d}^k \cdot p=p \}=d_j.
\end{equation*}
We set
\[
V_{d_j}=\bigoplus_{p \in \Sigma_{d_j}} V(p).
\]

Then, as seen in the example above, 
$T|_{V_{d_j}}$ can be represented by the RSP $S^j$ and its $d_j-1$ copies generated by the $C_d$-action on it.
We denote the collection of $S^j$ and its $d_j-1$ copies by $S^j_{d_j}$.
If $d_j$ equals 1, we simply write $S^j_{d_j}$ as $S^j$.

In this way, the HTL form (\ref{eq:direct_summand}) can be expressed as $S^1_{d_1} \ldots S^r_{d_r}$, 
which we call the \textit{spectral type of the singular point}.
The tuple of the spectral types of all singular points is called the \textit{spectral type of the equation}.
\begin{eg}
We provide more examples.
The spectral types of
\begin{equation*}
{\small
\begin{matrix}
x=0 \ \left( \frac12 \right) \\
\overbrace{
\begin{array}{cc}
a & \alpha\\
a & \alpha \\
-a & \alpha \\
-a & \alpha
\end{array}}
\end{matrix},\quad
\begin{matrix}
x=0 \ \left( \frac12 \right) \\
\overbrace{
\begin{array}{cc}
a & \alpha\\
-a & \alpha \\
b & \beta \\
-b & \beta
\end{array}}
\end{matrix},\quad
\begin{matrix}
x=0 \ \left( \frac12 \right) \\
\overbrace{
\begin{array}{cc}
a & \alpha\\
-a & \alpha \\
0 & \beta \\
0 &  \gamma
\end{array}}
\end{matrix},\quad
\begin{matrix}
x=0 \ \left( \frac13 \right) \\
\overbrace{
\begin{array}{cc}
a               & \alpha \\
\omega a    & \alpha\\
\omega^2 a & \alpha \\
0 & \beta
\end{array}}
\end{matrix},\quad
\begin{matrix}
  x=0 \ \left( \frac14 \right) \\
\overbrace{
\begin{array}{cc}
a              & \alpha \\
\sqrt{-1}a  & \alpha \\
-a             & \alpha \\
-\sqrt{-1}a & \alpha
\end{array}}
\end{matrix}
}
\end{equation*}
(where $\omega$ is a cube root of unity)
are $(2)_2$, $(1)_2(1)_2$, $(1)_2 11$, $(1)_3 1$, and $(1)_4$ respectively.
\qed
\end{eg}

\subsection{Degeneration of HTL forms}
As was mentioned in Section~1,
linear differential equations admit two kinds of degeneration;
namely, confluence of singular points and degeneration of HTL forms.

In this subsection, we illustrate the degeneration of an HTL form
which is caused by the degeneration of a Jordan canonical form.

\begin{rem}
The degenerations of HTL forms treated in this paper
correspond to degenerations of Jordan canonical forms.
However this is not true in general, see \cite{K2, K3}.
\qed
\end{rem}

First, we see an example of the degeneration of a Jordan canonical form. The following matrix
\begin{equation}\label{eq:eg_mat}
\begin{pmatrix}
\eta & 1 \\
0 & 0
\end{pmatrix}
\end{equation}
is diagonalizable provided that $\eta \neq 0$, and
when $\eta=0$ the matrix is nilpotent.
In other words, the semisimple matrix (\ref{eq:eg_mat}) degenerates to a nilpotent matrix
when $\eta$ tends to 0.

A matrix $g$ that diagonalizes the matrix (\ref{eq:eg_mat}) is of the form
\begin{equation}\label{eq:diag_gauge}
g=
\begin{pmatrix}
g_1 & -g_2/\eta \\
0 & g_2
\end{pmatrix}
\quad (g_1, g_2 \in \mathbb{C}^\times),
\end{equation}
that is,
\[
\begin{pmatrix}
\eta & 1 \\
0 & 0
\end{pmatrix}=
g
\begin{pmatrix}
\eta & 0 \\
0 & 0
\end{pmatrix}
g^{-1}
\]
holds.
Although the freedom $g_1$ and $g_2$ is redundant in this case,
we have retained both for later consideration.

Now we consider the following system:
\begin{equation}
\frac{d Y}{d x}=
\left(
\frac{A_0^{(-1)}}{x^2}+\frac{A_0^{(0)}}{x}+A_\infty
\right)Y.
\end{equation}
Here
$A_0^{(-1)}$, $A_0^{(0)}$, and $A_\infty$ are given as follows:
\begin{align*}
A_0^{(-1)}&=
\begin{pmatrix}
I_2 \\
P
\end{pmatrix}
\begin{pmatrix}
t(1-P) & tI_2
\end{pmatrix},\quad
A_0^{(0)}=
\begin{pmatrix}
-\theta^\infty_1I_2 & -Q \\
-Z & -\Theta
\end{pmatrix},\quad
A_\infty=
\begin{pmatrix}
-I_2 & O \\
O & O
\end{pmatrix},\\
Z&=(QP+\theta^0+2\theta^{\infty}_1)P-(QP+\theta^0+\theta^{\infty}_1),\quad
\Theta=
\begin{pmatrix}
\theta^\infty_2 & \\
 & \theta^\infty_3
\end{pmatrix}.
\end{align*}
Note that $Q$ and $P$ are assumed to satisfy the relation $[P, Q]=\theta^0+\theta^\infty_1+\Theta$.

The Riemann scheme of the system is given by
\[
\left(
\begin{array}{cc}
  x=0 & x=\infty \\
\overbrace{\begin{array}{cc}
     0 & 0 \\
     0 & 0 \\
     t & \theta^0 \\
     t & \theta^0
           \end{array}}
& 
\overbrace{\begin{array}{cc}
     1 & \theta^\infty_1 \\
     1 & \theta^\infty_1 \\
     0 & \theta^\infty_2 \\
     0 & \theta^\infty_3
           \end{array}} 
\end{array}
\right)
\]
and thus the spectral type of which is $(2)(2),(2)(11)$ (\cite{KNS}).

By changing the independent variable as $x=\varepsilon \tilde{x}$, we have
\begin{equation}
\frac{d Y}{d \tilde{x}}=
\left(
\frac{\varepsilon^{-1}A_0^{(-1)}}{\tilde{x}^2}+\frac{A_0^{(0)}}{\tilde{x}}+\varepsilon A_{\infty}
\right)Y.
\end{equation}

In a similar manner to (\ref{eq:diag_gauge}), we put
\begin{equation*}
G=
\begin{pmatrix}
G_1 & \varepsilon^{-1}G_2 \\
O & G_2
\end{pmatrix}
\end{equation*}
where $G_1$, $G_2$ are $2\times 2$ matrices.
The arbitrarity such as $G_1$ and $G_2$ will be used to eliminate negative powers of $\varepsilon$ or
simplify the system which result from the limit $\varepsilon \to 0$.
Using $G$ we have
\begin{equation}\label{eq:oshima-form}
G(\varepsilon A_\infty)G^{-1}=
\begin{pmatrix}
-\varepsilon I & I \\
O & O
\end{pmatrix},
\end{equation}
which tends to a nilpotent matrix as $\varepsilon \to 0$.
\begin{rem}
The matrices (\ref{eq:eg_mat}) and (the right-hand side of) (\ref{eq:oshima-form}) are
normal forms of matrices introduced by Oshima \cite{Os2005}.
\end{rem}

By changing the dependent variable as $Y=G^{-1}\tilde{Y}$,
we have
\begin{equation}
\frac{d \tilde{Y}}{d \tilde{x}}=G
\left(
\frac{\varepsilon^{-1}A_0^{(-1)}}{\tilde{x}^2}+\frac{A_0^{(0)}}{\tilde{x}}+\varepsilon A_{\infty}
\right)G^{-1}\tilde{Y}.
\end{equation}
Here we introduce the following transformations:
\begin{equation}
\begin{array}{l}
\theta^\infty_1=\varepsilon^{-1},\ \theta^\infty_2=\tilde{\theta}^\infty_2-\varepsilon^{-1},\ 
\theta^\infty_3=\tilde{\theta}^\infty_3-\varepsilon^{-1},\ t=\varepsilon\tilde{t},\\
Q=\varepsilon \tilde{Q}-(\theta^0\varepsilon+1)\tilde{P}^{-1},\ 
P=\varepsilon^{-1}\tilde{P}.
\end{array}
\end{equation}
Moreover, we choose $G_1=-\varepsilon P$, $G_2=I$.
Then, by a direct calculation, we have
\begin{align*}
G(\varepsilon^{-1}A_0^{(-1)})G^{-1}&=
\tilde{t}
\begin{pmatrix}
O \\
I
\end{pmatrix}
\begin{pmatrix}
\tilde{P} & I
\end{pmatrix}
+O(\varepsilon), \\
GA_0^{(0)}G^{-1}&=
\begin{pmatrix}
\tilde{Q}\tilde{P} & \tilde{Q} \\
I & -\tilde{P}\tilde{Q}+\theta^0
\end{pmatrix}
+O(\varepsilon).
\end{align*}
Taking the limit $\varepsilon \to 0$, we obtain a system of linear differential equations
which has a ramified singularity at $x=\infty$.
In fact, this is the system (\ref{eq:eg_shear}), which we have treated in Section~\ref{sec:shearing}.
This is an explicit description of the degeneration of an HTL form.
In terms of spectral types, it is expressed as $(2)(11), (2)(2) \to (11)_2, (2)(2)$.



\section{Lax pairs of degenerate matrix Painlev\'e systems}\label{sec:Lax_pairs}
The sixth matrix Painlev\'e system is derived from the isomonodromic deformation of the following Fuchsian system:
\begin{equation}\label{eq:Fuchs_mp}
\frac{dY}{dx}=
\left(
\frac{A_0}{x}+\frac{A_1}{x-1}+\frac{A_t}{x-t}
\right)Y
\end{equation}
where $A_0$, $A_1$, and $A_t$ are $4 \times 4$ matrices satisfying the following conditions
\begin{equation}\label{eq:eigen_condition}
A_0 \sim
\begin{pmatrix}
O_2 & O_2 \\
O_2 & \theta^0 I_2
\end{pmatrix},\ 
A_1 \sim
\begin{pmatrix}
O_2 & O_2 \\
O_2 & \theta^1 I_2
\end{pmatrix},\ 
A_t \sim
\begin{pmatrix}
O_2 & O_2 \\
O_2 & \theta^t I_2
\end{pmatrix},
\end{equation}
and
\begin{equation}\label{residue_infty}
A_\infty:=-(A_0+A_1+A_t)=
\mathrm{diag}(\theta^{\infty}_1,\theta^{\infty}_1,\theta^{\infty}_2,\theta^{\infty}_3).
\end{equation}
Thus the spectral type of the Fuchsian system (\ref{eq:Fuchs_mp}) is $22,22,22,211$.
Taking the trace of (\ref{residue_infty}), we have the Fuchs relation
\[
2(\theta^0+\theta^1+\theta^t+\theta^\infty_1)
+\theta^\infty_2+\theta^\infty_3=0.
\]
The explicit parametrization of (\ref{eq:Fuchs_mp}) is as follows~(\cite{K, KNS}):
\begin{equation}
\begin{split}
A_{\xi}&=
(U \oplus I_2)^{-1}X^{-1}\hat{A}_{\xi}X(U \oplus I_2)
\quad(\xi=0,1,t),\\
\hat{A}_0&=
\begin{pmatrix}
I_2 \\
O
\end{pmatrix}
\begin{pmatrix}
\theta^0I_2 & \frac{1}{t}Q-I_2
\end{pmatrix},\quad
\hat{A}_1=
\begin{pmatrix}
I_2 \\
PQ-\Theta
\end{pmatrix}
\begin{pmatrix}
\theta^1I_2-PQ+\Theta & I_2
\end{pmatrix},\\
\hat{A}_t&=
\begin{pmatrix}
I_2 \\
tP
\end{pmatrix}
\begin{pmatrix}
\theta^tI_2+QP & -\frac{1}{t}Q
\end{pmatrix}, \quad U \in \mathrm{GL}(2),
\end{split}
\end{equation}
where the matrices $Q$, $P$, and $\Theta$ are given by
\begin{align}\label{eq:canonical_var}
Q=
\begin{pmatrix}
q_1 & u\\
-q_2/u & q_1
\end{pmatrix},
\quad
P=
\begin{pmatrix}
p_1/2 & -p_2u\\
(p_2 q_2-\theta-\theta^\infty_1-\theta^\infty_2)/u & p_1/2
\end{pmatrix},
\quad
\Theta=
\begin{pmatrix}
\theta^{\infty}_2 & \\
 & \theta^{\infty}_3
\end{pmatrix}.
\end{align}
Here we have set $\theta=\theta^0+\theta^1+\theta^t$.
Note that $P$ and $Q$ satisfies $[P, Q]=(\theta+\theta^\infty_1)I_2+\Theta$.
The matrix $X$ is given by
$X=
\begin{pmatrix}
I_2 & O \\
Z & I_2
\end{pmatrix}$
where
\begin{align*}
Z&=(\theta^{\infty}_1-\Theta)^{-1}
[-\theta^1(QP+\theta+\theta^{\infty}_1)
+(QP+\theta+\theta^{\infty}_1)^2
-t(PQ+\theta^t)P].
\end{align*}
Note that the gauge parametrization is slightly different from that in \cite{KNS}, see also Appendix~\ref{sec:appendix_data} of the present paper.

As mentioned in Section~\ref{sec:intro},
the isomonodromic deformation equation of (\ref{eq:Fuchs_mp}) is equivalent to the compatibility condition of the following Lax pair:
\begin{equation}\label{eq:Lax_mp6}
\left\{
\begin{aligned}
\frac{\partial Y}{\partial x}&=
\left(
\frac{A_0}{x}+\frac{A_1}{x-1}+\frac{A_t}{x-t}
\right)Y,\\
\frac{\partial Y}{\partial t}&=-\frac{A_t}{x-t}Y.
\end{aligned}
\right.
\end{equation}



Moreover, the compatibility condition of the Lax pair (\ref{eq:Lax_mp6}) is equivalent to
{\small
\begin{align}
t(t-1)\frac{dQ}{dt}&=(Q-t)PQ(Q-1)+Q(Q-1)P(Q-t) \nonumber \\
&\quad+(\theta^0+1)Q(Q-1)+(\theta+2\theta^\infty_1-1)Q(Q-t)+\theta^t(Q-1)(Q-t), \label{eq:MP6-q}\\
t(t-1)\frac{dP}{dt}&=-(Q-1)P(Q-t)P-P(Q-t)PQ-PQ(Q-1)P \nonumber \\
&\quad-\left[(\theta^0+1)\{P(Q-1)+QP\}+(\theta+2\theta^\infty_1-1)\{P(Q-t)+QP\}+\theta^t\{P(Q-t)+(Q-1)P\}\right] \nonumber \\
&\quad-(\theta+\theta^\infty_1)(\theta^0+\theta^t+\theta^\infty_1), \label{eq:MP6-p}\\
t(t-1)\frac{dU}{dt}&=\left\{-\theta^1Q+(Q-t)(PQ+QP)+2(\theta+\theta^\infty_1)Q-\theta^t t\right\}U.
\end{align}
}
Equations (\ref{eq:MP6-q}) and (\ref{eq:MP6-p}) can be written in the following form:
\begin{align}
\frac{dq_i}{dt}&=\frac{\partial H^{\mathrm{Mat}}_{\mathrm{VI}}}{\partial p_i},\quad
\frac{dp_i}{dt}=-\frac{\partial H^{\mathrm{Mat}}_{\mathrm{VI}}}{\partial q_i}\quad (i=1,2),\\
\frac{du}{dt}&=-2q_1(q_1-1)(q_1-t)p_2
+\{q_1(q_1-1)+q_1(q_1-t)+(q_1-1)(q_1-t)-q_2\}p_1\nonumber\\
&\quad+(2p_2q_2-\theta^1-\theta^t-2\theta^\infty_2)q_1
+(2p_2q_2-\theta^0-2\theta^1-\theta^t-2\theta^\infty_1-2\theta^\infty_2+1)(q_1-1)\nonumber\\
&\quad+(2p_2q_2+\theta^0+\theta^1+2\theta^t+2\theta^\infty_1-1)(q_1-t)
\end{align}
where the Hamiltonian is given by
\begin{multline}
t(t-1)H^{\mathrm{Mat}}_{\mathrm{VI}}\left({-\theta^0-\theta^t-\theta^\infty_1,-\theta^1,\theta^t \atop
\theta^0+1,\theta+\theta^\infty_1+\theta^\infty_2};t;
{q_1, p_1\atop q_2, p_2}\right)\\=
\mathrm{tr}\Big[Q(Q-1)(Q-t)P^2+
\{(\theta^0+1 -(\theta+\theta^\infty_1+\Theta))Q(Q-1)+\theta^t(Q-1)(Q-t)
\\
+(\theta+2\theta^\infty_1-1)Q(Q-t)\}P
+(\theta+\theta^\infty_1)(\theta^0+\theta^t+\theta^\infty_1)Q\Big].
\end{multline}
We call (\ref{eq:MP6-q}) and (\ref{eq:MP6-p}) the non-abelian description of the sixth matrix Painlev\'e system
(see Appendix~\ref{sec:non-abel}).

From~(\ref{eq:canonical_var}), we can write the symplectic form as $\mathrm{tr}(dP \wedge dQ)$.

In this section, we list the linear systems of ramified type which are degenerated from~(\ref{eq:Fuchs_mp})
together with their deformations (i.e. systems of $t$-direction) and associated Hamiltonians.

We provide in advance the Hamilonians of the matrix Painlev\'e systems:
{\allowdisplaybreaks
\begin{align}
&t(t-1)H^{\mathrm{Mat}}_{\mathrm{VI}}
\left({\alpha, \beta, \gamma \atop \delta, \zeta}; t; Q,P\right) \nonumber \\
&=\mathrm{tr}[Q(Q-1)(Q-t)P^2 \nonumber \\
&\quad+\{(\delta-\zeta K)Q(Q-1)-(2\alpha+\beta+\gamma+\delta)Q(Q-t)
+\gamma(Q-1)(Q-t)\}P \nonumber \\
&\quad+\alpha(\alpha+\beta)Q],\\
&tH^{\mathrm{Mat}}_{\mathrm{V}}
\left({\alpha, \beta \atop \gamma, \zeta}; t; Q,P\right)=
\mathrm{tr}[P(P+t)Q(Q-1)+\beta PQ+\gamma P-(\alpha+\gamma) tQ],\\
&H^{\mathrm{Mat}}_{\mathrm{IV}}(\alpha, \beta, \zeta;t;Q,P)=\mathrm{tr}[PQ(P-Q-t)+\beta P+\alpha Q],\\
&tH^{\mathrm{Mat}}_{\mathrm{III}(D_6)}(\alpha, \beta, \zeta;t;Q,P)=\mathrm{tr}[P^2Q^2-(Q^2-\beta Q-t)P-\alpha Q],\\
&tH^{\mathrm{Mat}}_{\mathrm{III}(D_7)}(\alpha, \zeta;t;Q,P)=\mathrm{tr}[P^2Q^2+\alpha PQ+tP+Q],\\
&tH^{\mathrm{Mat}}_{\mathrm{III}(D_8)}(\zeta;t;Q,P)=\mathrm{tr}[P^2Q^2+PQ-Q-tQ^{-1}],\\
&H^{\mathrm{Mat}}_{\mathrm{II}}(\alpha, \zeta;t;Q,P)=\mathrm{tr}[P^2-(Q^2+t)P-\alpha Q],\\
&H^{\mathrm{Mat}}_{\mathrm{I}}(\zeta;t;Q,P)=\mathrm{tr}[P^2-Q^3-tQ].
\end{align}
}
Here the parameter $\zeta$
is included in such a manner as $[P, Q]=\zeta K$ where $K=\mathrm{diag}(1,-1)$.

\subsection{Singularity pattern $\frac32+1+1$}
\subsubsection{Spectral type $(2)_2, 22, 211$}
The Riemann scheme is given by
\[
\left(
\begin{array}{ccc}
  x=0 & x=1\ \left( \frac12 \right) & x=\infty \\
\begin{array}{c} 0 \\ 0 \\ \theta^0 \\ \theta^0 \end{array} &
\overbrace{\begin{array}{cc}
\sqrt{t} & 0\\
\sqrt{t} & 0\\
-\sqrt{t} & 0\\
-\sqrt{t} & 0
      \end{array}}&
\begin{array}{c} \theta^{\infty}_1 \\ \theta^{\infty}_1 \\ \theta^{\infty}_2 \\ \theta^{\infty}_3 \end{array} 
\end{array}
\right) ,
\]
and the Fuchs-Hukuhara relation is written as
$2\theta^0+2\theta_1^\infty+\theta_2^\infty+\theta_3^\infty=0$.

The Lax pair is expressed as
\begin{equation}\label{eq:Lax(2)_2,22,211}
\left\{
\begin{aligned}
\frac{\partial Y}{\partial x}&=
\left(
\frac{A_0^{(0)}}{x}+\frac{A_1^{(1)}}{(x-1)^2}+\frac{A_1^{(0)}}{x-1}
\right)Y ,\\
\frac{\partial Y}{\partial t}&=\left(\frac{-\frac{1}{t}A_1^{(1)}}{x-1}\right)Y .
\end{aligned}
\right.
\end{equation}
Here
$A_0^{(0)}$, $A_1^{(1)}$, and $A_1^{(0)}$ are given as follows:
\begin{align*}
A_{\xi}^{(k)}&=
(U \oplus I_2)^{-1}\hat{A}_{\xi}^{(k)}(U \oplus I_2),\\
\hat{A}_1^{(1)}&=
G_1
\begin{pmatrix}
O & -tI \\
O & O
\end{pmatrix}
G_1^{-1}
=
\begin{pmatrix}
I_2 \\
-\frac{1}{t}Z
\end{pmatrix}
\begin{pmatrix}
-Z & -tI 
\end{pmatrix},\\
\hat{A}_1^{(0)}&=
G_1
\begin{pmatrix}
PQ-\theta^\infty_1 & tP\\
I_2 & -PQ+\theta^\infty_1
\end{pmatrix}
G_1^{-1},\\
\hat{A}_0^{(0)}&=
\begin{pmatrix}
P \\
-\frac{1}{t}(ZP+QP+\theta^0)
\end{pmatrix}
\begin{pmatrix}
-Z-Q & -tI
\end{pmatrix},\\
Z&=(\theta^{\infty}_1-\Theta)^{-1}
( -QPQ-\theta^0 Q-t ),\quad
G_1=
\begin{pmatrix}
I_2 & O_2 \\
-\frac{1}{t}Z & I_2
\end{pmatrix}.
\end{align*}
Here $Q$, $P$, and $\Theta$ are
\begin{align*}
Q=
\begin{pmatrix}
q_1 & u \\
-q_2/u & q_1
\end{pmatrix},\quad
P=
\begin{pmatrix}
p_1/2 & -p_2 u\\
(p_2 q_2-\theta^0-\theta^\infty_1-\theta^\infty_2)/u & p_1/2
\end{pmatrix},\quad
\Theta=
\begin{pmatrix}
\theta^{\infty}_2 & \\
 & \theta^{\infty}_3
\end{pmatrix}.
\end{align*}

The compatibility condition of (\ref{eq:Lax(2)_2,22,211}) is equivalent to
{\small
\begin{align}
t\frac{dQ}{dt}&=2QPQ-Q^2-(2\theta^\infty_1-1)Q+t, \label{eq:(2)_2,22,211q}\\
t\frac{dP}{dt}&=-2PQP+PQ+QP+(2\theta^\infty_1-1)P+\theta^0,  \label{eq:(2)_2,22,211p}\\
t\frac{dU}{dt}U^{-1}&=-2PQ+2\theta^\infty_1.
\end{align}
}
Equations (\ref{eq:(2)_2,22,211q}) and (\ref{eq:(2)_2,22,211p}) can be written in the following form:
\begin{align}
&\frac{dq_i}{dt}=\frac{\partial H^{\mathrm{Mat}}_{\mathrm{III}(D_6)}}{\partial p_i},\quad
\frac{dp_i}{dt}=-\frac{\partial H^{\mathrm{Mat}}_{\mathrm{III}(D_6)}}{\partial q_i}\quad (i=1,2), \\
&t\frac{1}{u}\frac{du}{dt}=-2q_1(q_1p_2+1)+2(p_1q_1+p_2q_2)-2(\theta^0+2\theta^\infty_1+\theta^\infty_2)+1
\end{align}
where the Hamiltonian is given by
\begin{align}
&tH^{\mathrm{Mat}}_{\mathrm{III}(D_6)}\left({\theta^0, -2\theta^\infty_1+1, 
\theta^0+\theta^\infty_1+\theta^\infty_2};t;{q_1, p_1 \atop q_2, p_2}\right)\\
&=\mathrm{tr}\Big[ P^2Q^2-\{Q^2+(2\theta^\infty_1-1)Q-t\}P-\theta^0Q \Big].\nonumber
\end{align}

\subsubsection{Spectral type $(11)_2, 22, 22$}
The Riemann scheme is given by
\[
\left(
\begin{array}{ccc}
  x=0 & x=1 & x=\infty \ \left( \frac12 \right)\\
\begin{array}{c} 0 \\ 0 \\ \theta^0 \\ \theta^0 \end{array} &
\begin{array}{c} 0 \\ 0 \\ \theta^1 \\ \theta^1 \end{array} &
\overbrace{\begin{array}{cc}
     \sqrt{t}   &   \theta^\infty_2/2 \\
     \sqrt{t}   &   \theta^\infty_3/2 \\
     -\sqrt{t} & \theta^\infty_2/2 \\
     -\sqrt{t} & \theta^\infty_3/2
           \end{array}}
\end{array}
\right) ,
\]
and the Fuchs-Hukuhara relation is written as
$2\theta^0 +2\theta^1+\theta_2^\infty+\theta_3^\infty =0$.

The Lax pair is expressed as
\begin{equation}\label{eq:Lax(11)_2,22,22}
\left\{
\begin{aligned}
\frac{\partial Y}{\partial x}&=
\left(A_\infty+\frac{A_0}{x}+\frac{A_1}{x-1}
\right)Y, \\
\frac{\partial Y}{\partial t}&=(N x+B_1)Y ,
\end{aligned}
\right.
\end{equation}
where
\begin{align*}
A_\infty&=tN,\quad
N=
\begin{pmatrix}
O_2 & I_2 \\
O_2 & O_2
\end{pmatrix},\quad
A_0=
\begin{pmatrix}
O_2 \\
I_2
\end{pmatrix}
\begin{pmatrix}
I_2-P & \theta^0 I_2
\end{pmatrix},\quad
A_1=
\begin{pmatrix}
QP+\theta^1\\
P
\end{pmatrix}
\begin{pmatrix}
I_2 & -Q
\end{pmatrix}.
\end{align*}
Furthermore the matrix $B_1$ is given by
\begin{align*}
B_1&=\frac1t
\begin{pmatrix}
PQ-\theta^0 & O_2\\
I_2 & -QP-\theta^1
\end{pmatrix}.
\end{align*}
Here $Q$, $P$, and $\Theta$ are
\begin{align*}
Q=
\begin{pmatrix}
q_1 & u \\
-q_2/u & q_1
\end{pmatrix},\quad
P=
\begin{pmatrix}
p_1/2 & -p_2 u\\
(p_2 q_2-\theta^0-\theta^1-\theta^\infty_2)/u & p_1/2
\end{pmatrix},\quad
\Theta=
\begin{pmatrix}
\theta^{\infty}_2 & \\
 & \theta^{\infty}_3
\end{pmatrix}.
\end{align*}

The compatibility condition of (\ref{eq:Lax(11)_2,22,22}) is equivalent to
\begin{align}
t\frac{dQ}{dt}&=PQ^2+Q^2P-Q^2-(\theta^0-\theta^1)Q+t, \label{eq:(11)_2,22,22q}\\
t\frac{dP}{dt}&=-P^2Q-QP^2+PQ+QP+(\theta^0-\theta^1)P+\theta^1. \label{eq:(11)_2,22,22p}
\end{align}
These are apparently different from the expressions in Appendix~\ref{sec:non-abel}.
However, using the relation $[P, Q]=\theta^0+\theta^1+\Theta$ and a gauge transformation
\[
Q=G^{-1}\tilde{Q}G, \quad P=G^{-1}\tilde{P}G, \quad G=\mathrm{diag}(t^{-\theta^\infty_2}, t^{-\theta^\infty_3}),
\]
we can see that (\ref{eq:(11)_2,22,22q}) and (\ref{eq:(11)_2,22,22p}) are essentially the same as (\ref{eq:non_abel_III(D6)}) in Appendix~\ref{sec:non-abel}.

Equations (\ref{eq:(11)_2,22,22q}), (\ref{eq:(11)_2,22,22p}) can be written in the following form:
\begin{align}
&\frac{dq_i}{dt}=\frac{\partial H^{\mathrm{Mat}}_{\mathrm{III}(D_6)}}{\partial p_i},\quad
\frac{dp_i}{dt}=-\frac{\partial H^{\mathrm{Mat}}_{\mathrm{III}(D_6)}}{\partial q_i}\quad (i=1,2),\\
&t\frac{1}{u}\frac{du}{dt}=2(p_1q_1+p_2q_2)-2q_1(q_1p_2+1)-\theta^0+\theta^1.
\end{align}
The Hamiltonian is given by
\begin{equation}
tH^{\mathrm{Mat}}_{\mathrm{III}(D_6)}\left({\theta^1,\theta^1-\theta^0,\theta^0+\theta^1+\theta^\infty_2};t;
{q_1, p_1 \atop q_2, p_2}\right)=
\mathrm{tr}\Big[ P^2Q^2-(Q^2+(\theta^0-\theta^1)Q-t)P-\theta^1Q \Big].
\end{equation}

\subsection{Singularity pattern $\frac52+1$}
\subsubsection{Spectral type $(((2)))_2, 211$}\label{sec:(((2)))_2,211}
The Riemann scheme is given by
\[
\left(
\begin{array}{cc}
  x=0 \ \left( \frac12 \right) & x=\infty \\
\overbrace{\begin{array}{cccc}
   1& 0 & -t/2  & 0 \\
   1& 0 & -t/2 & 0 \\
  -1& 0 & t/2 & 0 \\
  -1& 0 & t/2 & 0
        \end{array}}
& \begin{array}{c} \theta^{\infty}_1 \\ \theta^{\infty}_1 \\ \theta^{\infty}_2 \\ \theta^{\infty}_3 \end{array}
\end{array}
\right) ,
\]
and the Fuchs-Hukuhara relation is written as
$2\theta_1^\infty+\theta_2^\infty+\theta_3^\infty=0$.

The Lax pair is expressed as
\begin{equation}\label{eq:Lax(((2)))_2,211}
\left\{
\begin{aligned}
\frac{\partial Y}{\partial x}&=
\left(
\frac{A_{2}}{x^3}+\frac{A_{1}}{x^2}+\frac{A_{0}}{x}
\right)Y ,\\
\frac{\partial Y}{\partial t}&=\frac{A_{2}}{x}Y .
\end{aligned}
\right.
\end{equation}
The matrices $A_0$, $A_1$, and $A_2$ are given as follows:
\begin{align*}
A_{k}&=
(U \oplus I_2)^{-1}\hat{A}_{k}(U \oplus I_2),\\
\hat{A}_{2}&=
G_0
\begin{pmatrix}
O & I_2 \\
O & O
\end{pmatrix}
G_0^{-1},\quad
\hat{A}_{1}=
G_0
\begin{pmatrix}
Q & -P \\
I & -Q
\end{pmatrix}
G_0^{-1},\quad
\hat{A}_{0}=-
\begin{pmatrix}
\theta^{\infty}_1I_2 & O \\
O & \Theta
\end{pmatrix},\\
G_0&=
\begin{pmatrix}
I & O \\
Z & I
\end{pmatrix},
\quad
Z=(\theta^{\infty}_1-\Theta)^{-1}(P-Q^2-t).
\end{align*}
Here, $Q$, $P$, and $\Theta$ are
\begin{align*}
Q=
\begin{pmatrix}
q_1 & u\\
-q_2/u & q_1
\end{pmatrix},
\quad
P=
\begin{pmatrix}
p_1/2 & -p_2 u\\
(p_2q_2-\theta^\infty_1-\theta^\infty_2)/u & p_1/2
\end{pmatrix},\quad
\Theta=
\begin{pmatrix}
\theta^{\infty}_2 & \\
 & \theta^{\infty}_3
\end{pmatrix}.
\end{align*}

The compatibility condition of (\ref{eq:Lax(((2)))_2,211}) is
\begin{align}
\frac{dQ}{dt}&=2P-Q^2-t, \quad
\frac{dP}{dt}=PQ+QP-2\theta^\infty_1+1, \label{eq:(((2)))_2,211} \\
\frac{dU}{dt}U^{-1}&=2Q.
\end{align}
Equations (\ref{eq:(((2)))_2,211}) are written in the following form:
\begin{align}
\frac{dq_i}{dt}=\frac{\partial H^{\mathrm{Mat}}_{\mathrm{II}}}{\partial p_i},\quad
\frac{dp_i}{dt}=-\frac{\partial H^{\mathrm{Mat}}_{\mathrm{II}}}{\partial q_i}\quad (i=1,2),\quad
\frac{1}{u}\frac{du}{dt}=-2(q_1+p_2),
\end{align}
where the Hamiltonian is given by
\begin{equation}
H^{\mathrm{Mat}}_{\mathrm{II}}\left(-2\theta^\infty_1+1, \theta^\infty_1+\theta^\infty_2 ;t;{q_1, p_1\atop q_2, p_2}\right)
 =\mathrm{tr}\Big[ P^2-(Q^2+t)P+(2\theta^\infty_1-1)Q \Big].
\end{equation}

\subsubsection{Spectral type $(((11)))_2, 22$}\label{sec:(((11)))_2,22}
The Riemann scheme is given by
\[
\left(
\begin{array}{cc}
  x=0 & x=\infty \ \left( \frac12 \right) \\
\begin{array}{c} 0 \\ 0 \\ \theta^0 \\ \theta^0\end{array}
&\overbrace{\begin{array}{cccc}
  1 & 0 & -t/2 & \theta^\infty_2/2 \\
  1 & 0 & -t/2 & \theta^\infty_3/2 \\
 -1& 0 & t/2 & \theta^\infty_2/2 \\
 -1& 0 & t/2 & \theta^\infty_3/2
        	\end{array}}
\end{array}
\right) ,
\]
and the Fuchs-Hukuhara relation is written as
$2\theta^0 +\theta_2^\infty +\theta_3^\infty =0$.

The Lax pair is expressed as
\begin{equation}\label{eq:Lax(((11)))_2,22}
\left\{
\begin{aligned}
\frac{\partial Y}{\partial x}&=
\left(
A_0x+A_1+\frac{A_2}{x}
\right)Y ,\\
\frac{\partial Y}{\partial t}&=(-A_0x+B_1)Y ,
\end{aligned}
\right.
\end{equation}
where
\begin{align*}
A_0&=
\begin{pmatrix}
O_2 & I_2\\
O_2 & O_2
\end{pmatrix},\quad
A_1=
\begin{pmatrix}
O_2 & P-tI_2\\
I_2 & O_2
\end{pmatrix},\quad
A_2=
\begin{pmatrix}
-Q \\
I_2
\end{pmatrix}
\begin{pmatrix}
-P & -PQ+\theta^0I_2
\end{pmatrix},\\
B_1&=
\begin{pmatrix}
O_2 & -2P+tI_2\\
-I_2 & O_2
\end{pmatrix}.
\end{align*}
Here $Q$, $P$, and $\Theta$ are
\begin{align*}
Q=
\begin{pmatrix}
q_1 & u \\
-q_2/u & q_1
\end{pmatrix},\quad
P=
\begin{pmatrix}
p_1/2 & -p_2 u \\
(p_2q_2-\theta^0-\theta^\infty_2)/u & p_1/2
\end{pmatrix},\quad
\Theta=
\begin{pmatrix}
\theta^\infty_2 & \\
 & \theta^\infty_3
\end{pmatrix}.
\end{align*}

The compatibility condition of (\ref{eq:Lax(((11)))_2,22}) is equivalent to
\begin{align}
\frac{dQ}{dt}=2P-Q^2-t, \quad \frac{dP}{dt}=PQ+QP-\theta^0,
\end{align}
and these are written in the following form:
\begin{align}
\frac{dq_i}{dt}=\frac{\partial H^{\mathrm{Mat}}_{\mathrm{II}}}{\partial p_i}, \quad
\frac{dp_i}{dt}=-\frac{\partial H^{\mathrm{Mat}}_{\mathrm{II}}}{\partial q_i}, \quad
\frac{1}{u}\frac{du}{dt}=-2(q_1+p_2).
\end{align}
Here the Hamiltonian is given by
\begin{equation}
H^{\mathrm{Mat}}_{\mathrm{II}}\left(-\theta^0, \theta^0+\theta^\infty_2 ;t;{q_1, p_1\atop q_2, p_2}\right)
 =\mathrm{tr}\Big[ P^2-(Q^2+t)P+\theta^0 Q \Big].
\end{equation}


\subsection{Singularity pattern $\frac32+2$}
\subsubsection{Spectral type $(2)_2, (2)(11)$}
The Riemann scheme is given by
\[
\left(
\begin{array}{cc}
  x=0 \ \left( \frac12 \right)& x=\infty \\
\overbrace{\begin{array}{cc}
     \sqrt{-t}  & 0 \\
     \sqrt{-t}  & 0 \\
     -\sqrt{-t} & 0 \\
     -\sqrt{-t} & 0
           \end{array}}
& 
\overbrace{\begin{array}{cc}
     0 & \theta^\infty_1 \\
     0 & \theta^\infty_1 \\
    -1 & \theta^\infty_2 \\
    -1 & \theta^\infty_3
           \end{array}} 
\end{array}
\right) ,
\]
and the Fuchs-Hukuhara relation is written as
$2\theta_1^\infty+\theta_2^\infty+\theta_3^\infty=0$.

The Lax pair is expressed as
\begin{equation}\label{eq:Lax(2)_2,(2)(11)}
\left\{
\begin{aligned}
\frac{\partial Y}{\partial x}&=
\left(
\frac{A_0}{x^2}+\frac{A_1}{x}+A_2
\right)Y ,\\
\frac{\partial Y}{\partial
 t}&=-\frac{1}{x}\left(\frac{A_0}{t}\right)Y,
\end{aligned}
\right.
\end{equation}
where
$A_0$, $A_1$, and $A_2$ are given as follows:
\begin{align*}
A_{\xi}&=
(U \oplus I_2)^{-1}\hat{A}_{\xi}(U \oplus I_2),\\
\hat{A}_0&=t
\begin{pmatrix}
I_2 \\
P
\end{pmatrix}
\begin{pmatrix}
-P & I_2
\end{pmatrix},\quad
\hat{A}_1=
\begin{pmatrix}
-\theta^\infty_1I_2 & -Q \\
-Z & -\Theta
\end{pmatrix},\quad
\hat{A}_2=
\begin{pmatrix}
O & O \\
O & I_2
\end{pmatrix},\\
Z&=(QP+2\theta^{\infty}_1)P+I.
\end{align*}
Here $Q$, $P$, and $\Theta$ are
\begin{align*}
Q=
\begin{pmatrix}
q_1 & u\\
-q_2/u & q_1
\end{pmatrix},
\quad
P=
\begin{pmatrix}
p_1/2 & -p_2 u\\
(p_2q_2-\theta^\infty_1-\theta^\infty_2)/u & p_1/2
\end{pmatrix},\quad
\Theta=
\begin{pmatrix}
\theta^{\infty}_2 & \\
 & \theta^{\infty}_3
\end{pmatrix}.
\end{align*}

The compatibility condition of (\ref{eq:Lax(2)_2,(2)(11)}) is
\begin{align}
&t\frac{dQ}{dt}=2QPQ+2\theta^\infty_1Q+t, \quad
t\frac{dP}{dt}=-2PQP-2\theta^\infty_1P-1, \label{eq:(2)_2,(2)(11)qp}\\
&t\frac{dU}{dt}U^{-1}=2(QP+\theta^\infty_1).
\end{align}
Equations~(\ref{eq:(2)_2,(2)(11)qp}) are written in the following form:
\begin{align}
&\frac{dq_i}{dt}=\frac{\partial H^{\mathrm{Mat}}_{\mathrm{III}(D_7)}}{\partial p_i},\quad
\frac{dp_i}{dt}=-\frac{\partial H^{\mathrm{Mat}}_{\mathrm{III}(D_7)}}{\partial q_i}\quad (i=1,2),\\
&t\frac{1}{u}\frac{du}{dt}=2(p_1q_1+p_2q_2-p_2q_1^2-\theta^\infty_2).
\end{align}
The Hamiltonian is given by
\begin{align}
&tH^{\mathrm{Mat}}_{{\mathrm{III}}(D_7)}
\left({2\theta^\infty_1,\theta^\infty_1+\theta^\infty_2};t;
  {q_1,p_1\atop q_2,p_2}\right)
=\mathrm{tr}\Big[ P^2Q^2+2\theta^\infty_1PQ+tP+Q \Big].
\end{align}

\subsubsection{Spectral type $(11)_2, (2)(2)$}\label{sec:(11)_2,(2)(2)}
The Riemann scheme is given by
\[
\left(
\begin{array}{cc}
  x=0 & x=\infty \ \left( \frac12 \right)\\
\overbrace{\begin{array}{cc}
     0 & 0 \\
     0 & 0 \\
     t & \theta^0 \\
     t & \theta^0
           \end{array}}
& 
\overbrace{\begin{array}{cc}
     1 & \theta^\infty_2/2 \\
     1 & \theta^\infty_3/2 \\
     -1 & \theta^\infty_2/2 \\
     -1 & \theta^\infty_3/2
           \end{array}} 
\end{array}
\right) ,
\]
and the Fuchs-Hukuhara relation is written as
$2\theta^0+\theta_2^\infty+\theta_3^\infty=0$.

The Lax pair is expressed as
\begin{equation}\label{eq:Lax(11)_2,(2)(2)}
\left\{
\begin{aligned}
\frac{\partial Y}{\partial x}&=\left( \frac{A_0}{x^2}+\frac{A_1}{x}+A_2 \right)Y, \\
\frac{\partial Y}{\partial t}&=\left(B_0+\frac{B_1}{x}\right)Y,
\end{aligned}
\right.
\end{equation}
where $A_0$, $A_1$, $A_2$, $B_0$, and $B_1$ are given as follows:
\begin{align*}
A_0&=t
\begin{pmatrix}
O_2 \\
I_2
\end{pmatrix}
\begin{pmatrix}
P &  I_2
\end{pmatrix},\quad
A_1=
\begin{pmatrix}
QP & Q \\
I_2 & -PQ+\theta^0
\end{pmatrix},\quad
A_2=
\begin{pmatrix}
O & I \\
O & O
\end{pmatrix},\\
B_0&=
-\frac1t
\begin{pmatrix}
O & Q \\
O & O
\end{pmatrix},\quad
B_1=
-\begin{pmatrix}
O_2 \\
I_2
\end{pmatrix}
\begin{pmatrix}
P & I_2
\end{pmatrix}.
\end{align*}
Here $Q$, $P$, and $\Theta$ are
\begin{align*}
Q=
\begin{pmatrix}
q_1 & u\\
-q_2/u & q_1
\end{pmatrix},
\quad
P=
\begin{pmatrix}
p_1/2 & -p_2 u\\
(p_2q_2-\theta^0-\theta^\infty_2)/u & p_1/2
\end{pmatrix},\quad
\Theta=
\begin{pmatrix}
\theta^{\infty}_2 & \\
 & \theta^{\infty}_3
\end{pmatrix}.
\end{align*}

The compatibility condition of (\ref{eq:Lax(11)_2,(2)(2)}) is
\begin{align}
t\frac{dQ}{dt}=2QPQ-\theta^0 Q+t,\quad t\frac{dP}{dt}=-2PQP+\theta^0 P-1,
\end{align}
and these are written in the following form:
\begin{align}
\frac{dq_i}{dt}=\frac{\partial H^{\mathrm{Mat}}_{\mathrm{III}(D_7)}}{\partial p_i}, \quad
\frac{dp_i}{dt}=-\frac{\partial H^{\mathrm{Mat}}_{\mathrm{III}(D_7)}}{\partial q_i}, \quad
t\frac{1}{u}\frac{du}{dt}=2p_1q_1+2p_2q_2-2p_2q_1^2-3\theta^0-2\theta^\infty_2.
\end{align}
Here the Hamiltonian is given by
\begin{equation}
tH^{\mathrm{Mat}}_{\mathrm{III}(D_7)}\left(-\theta^0, \theta^0+\theta^\infty_2 ;t;{q_1, p_1\atop q_2, p_2}\right)
 =\mathrm{tr}\Big[ P^2Q^2-\theta^0 PQ+tP+Q \Big].
\end{equation}

\subsection{Singularity pattern $\frac72$}
\subsubsection{Spectral type $(((((11)))))_2$}
The Riemann scheme is given by
\[
\left(
\begin{array}{c}
 x=\infty \ \left( \frac12 \right) \\
\overbrace{\begin{array}{cccccc}
 1 & 0 & 0 & 0 & t/2 & \theta^\infty_2/2\\
 1 & 0 & 0 & 0 & t/2 & \theta^\infty_3/2\\
 -1 & 0 & 0 & 0 & -t/2 & \theta^\infty_2/2\\
 -1 & 0 & 0 & 0 & -t/2 & \theta^\infty_3/2
        	\end{array}}
\end{array}
\right) ,
\]
and the Fuchs-Hukuhara relation is written as
$\theta_2^\infty+\theta_3^\infty =0$.

The Lax pair is expressed as
\begin{equation}\label{eq:Lax_(((((11)))))_2}
\left\{
\begin{aligned}
\frac{\partial Y}{\partial x}&=
\left(
A_0x^2+A_1x+A_2
\right)Y ,\\
\frac{\partial Y}{\partial t}&=(A_0x+B_1)Y ,
\end{aligned}
\right.
\end{equation}
where
\begin{align*}
A_0&=
\begin{pmatrix}
O & I \\
O & O
\end{pmatrix},\quad
A_1=
\begin{pmatrix}
O    & Q \\
I & O  
\end{pmatrix},\quad
A_2=
\begin{pmatrix}
-P & Q^2+t \\
-Q & P
\end{pmatrix},\\
B_1&=
\begin{pmatrix}
O & 2Q \\
I & O
\end{pmatrix}.
\end{align*}
Here $Q$, $P$, and $\Theta$ are
\begin{align*}
Q=
\begin{pmatrix}
q_1 & u \\
-q_2/u & q_1
\end{pmatrix},\quad
P=
\begin{pmatrix}
p_1/2 & -p_2 u \\
(p_2q_2-\theta^\infty_2)/u & p_1/2
\end{pmatrix},\quad
\Theta=
\begin{pmatrix}
\theta^\infty_2 &  \\
 & \theta^\infty_3
\end{pmatrix}.
\end{align*}

The compatibility condition of (\ref{eq:Lax_(((((11)))))_2}) is
\begin{equation}
\frac{dQ}{dt}=2P, \quad \frac{dP}{dt}=3Q^2+t,
\end{equation}
and they are written in the following form:
\begin{align}
\frac{dq_i}{dt}=\frac{\partial H^{\mathrm{Mat}}_{\mathrm{I}}}{\partial p_i}, \quad
\frac{dp_i}{dt}=-\frac{\partial H^{\mathrm{Mat}}_{\mathrm{I}}}{\partial q_i}, \quad
\frac{1}{u}\frac{du}{dt}=-2p_2.
\end{align}
The Hamiltonian is given by
\begin{equation}
H^{\mathrm{Mat}}_{\mathrm{I}}\left({\theta^\infty_2};t;{q_1,p_1\atop q_2,p_2}\right)=\mathrm{tr}\Big[ P^2-Q^3-tQ \Big].
\end{equation}

\subsection{Singularity pattern $\frac32+\frac32$}
\subsubsection{Spectral type $(2)_2, (11)_2$}
The Riemann scheme is given by
\[
\left(
\begin{array}{cc}
  x=0 \ \left( \frac12 \right) & x=\infty \ \left( \frac12 \right)\\
\overbrace{\begin{array}{cc}
\sqrt{t} & 0 \\
\sqrt{t} & 0 \\
-\sqrt{t} & 0 \\
-\sqrt{t} & 0
           \end{array}}
& 
\overbrace{\begin{array}{cc}
     1 & \theta^\infty_2/2 \\
     1 & \theta^\infty_3/2 \\
     -1 & \theta^\infty_2/2 \\
     -1 & \theta^\infty_3/2
           \end{array}} 
\end{array}
\right) ,
\]
and the Fuchs-Hukuhara relation is written as
$\theta_2^\infty+\theta_3^\infty=0$.

The Lax pair is expressed as
\begin{equation}\label{eq:Lax(2)_2,(11)_2}
\left\{
\begin{aligned}
\frac{\partial Y}{\partial x}&=\left( \frac{A_0}{x^2}+\frac{A_1}{x}+A_2 \right)Y, \\
\frac{\partial Y}{\partial t}&=\left(B_0+\frac{B_1}{x}\right)Y,
\end{aligned}
\right.
\end{equation}
where $A_0$, $A_1$, $A_2$, $B_0$, and $B_1$ are given as follows:
\begin{align*}
A_0&=
\begin{pmatrix}
O_2 & O_2 \\
-tQ^{-1} & O_2
\end{pmatrix},\quad
A_1=
\begin{pmatrix}
QP & -Q \\
I_2 & -PQ-I_2
\end{pmatrix},\quad
A_2=
\begin{pmatrix}
O & I \\
O & O
\end{pmatrix},\\
B_0&=
\frac1t
\begin{pmatrix}
O & Q \\
O & O
\end{pmatrix},\quad
B_1=
\begin{pmatrix}
O_2 & O_2 \\
Q^{-1} & O_2
\end{pmatrix}.
\end{align*}
Here $Q$, $P$, and $\Theta$ are
\begin{align*}
Q=
\begin{pmatrix}
q_1 & u\\
-q_2/u & q_1
\end{pmatrix},
\quad
P=
\begin{pmatrix}
p_1/2 & -p_2 u\\
(p_2q_2-\theta^\infty_2)/u & p_1/2
\end{pmatrix},\quad
\Theta=
\begin{pmatrix}
\theta^{\infty}_2 & \\
 & \theta^{\infty}_3
\end{pmatrix}.
\end{align*}

The compatibility condition of (\ref{eq:Lax(2)_2,(11)_2}) is
\begin{align}
t\frac{dQ}{dt}=2QPQ+Q,\quad t\frac{dP}{dt}=-2PQP-P+I-tQ^{-2},
\end{align}
and these are written in the following form:
\begin{align}
\frac{dq_i}{dt}=\frac{\partial H^{\mathrm{Mat}}_{\mathrm{III}(D_8)}}{\partial p_i}, \quad
\frac{dp_i}{dt}=-\frac{\partial H^{\mathrm{Mat}}_{\mathrm{III}(D_8)}}{\partial q_i}, \quad
t\frac{1}{u}\frac{du}{dt}=2p_1q_1+2p_2q_2-2p_2q_1^2-2\theta^\infty_2+1.
\end{align}
Here the Hamiltonian is given by
\begin{equation}
tH^{\mathrm{Mat}}_{\mathrm{III}(D_8)}\left( \theta^\infty_2 ;t;{q_1, p_1\atop q_2, p_2}\right)
 =\mathrm{tr}\Big[ P^2Q^2+PQ-Q-tQ^{-1} \Big].
\end{equation}

\section{Laplace transform}\label{sec:Laplace}
In this section, we describe correspondences of linear systems through the Laplace transform.
So far the configuration of the singular points of a linear system is not important.
However, when we consider the Laplace transform, the singular point $x=\infty$ is distinguished from
the other singular points.
Hence, only in this section, we indicate the spectral type corresponding to $x=\infty$ with $\infty$ over the spectral type.

There are five Hamiltonians in the degeneration scheme of the matrix Painlev\'e systems in Section~\ref{sec:intro} which have more than one
associated linear systems; that is, $H^{\mathrm{Mat}}_{\mathrm{V}}$, $H^{\mathrm{Mat}}_{\mathrm{IV}}$,
$H^{\mathrm{Mat}}_{\mathrm{III}(D_6)}$, $H^{\mathrm{Mat}}_{\mathrm{III}(D_7)}$,
and $H^{\mathrm{Mat}}_{\mathrm{II}}$.
The correspondences concerning $H^{\mathrm{Mat}}_{\mathrm{V}}$ and $H^{\mathrm{Mat}}_{\mathrm{IV}}$ are given in \cite{KNS}.
Here we see the correspondences concerning $H^{\mathrm{Mat}}_{\mathrm{III}(D_6)}$, $H^{\mathrm{Mat}}_{\mathrm{III}(D_7)}$,
and $H^{\mathrm{Mat}}_{\mathrm{II}}$.

In the case of linear systems of the following form:
\begin{equation}
\frac{d}{dx}Y=\left[ B\left(xI_l-T\right)^{-1}C+S\right] Y
\end{equation}
where $B$ is an $m \times l$ matrix and $C$ is an $l \times m$ matrix,
the correspondence is very symmetric~(\cite{B2, Hrd}).
To see this, we rewrite this equation as
\begin{equation}
\begin{pmatrix}
\frac{d}{dx}-S & -B \\
-C & xI_l-T
\end{pmatrix}
\begin{pmatrix}
  Y \\
  Z
\end{pmatrix}=0.
\end{equation}
Applying the Laplace transform $(x, d/dx)\mapsto (-d/d\xi , \xi)$, we have
\begin{equation}
\begin{pmatrix}
\xi -S & -B \\
-C & -\frac{d}{d\xi}-T
\end{pmatrix}
\begin{pmatrix}
  \hat{Y} \\
  \hat{Z}
\end{pmatrix}=0.
\end{equation}
Here we have expressed the transformed dependent variables using $\hat{~}$.
If we regard this equation as the equation of $\hat{Z}$, the equation reads
\begin{equation}
\frac{d}{d\xi}\hat{Z}=-\left[ C\left(\xi I_m-S\right)^{-1}B+T\right]\hat{Z}.
\end{equation}
Through this procedure, we see the following correspondences of spectral types:
\begin{align*}
H^{\mathrm{Mat}}_{\mathrm{III}(D_6)}&:\  \stackrel{\infty}{(2)}_2, 22, 211 \leftrightarrow (2)(2), \stackrel{\infty}{(2)(11)}, \quad
\stackrel{\infty}{(11)}_2, 22, 22 \leftrightarrow (2)(11), \stackrel{\infty}{(2)(2)}, \\
H^{\mathrm{Mat}}_{\mathrm{III}(D_7)}&:\  \stackrel{\infty}{(2)}_2,(2)(11) \leftrightarrow (2)(2), \stackrel{\infty}{(11)}_2.
\end{align*}

The correspondences concerning $H^{\mathrm{Mat}}_{\mathrm{II}}$ are rather complicated.
For an example, let us look at the system $(((11)))_2, 22$ (Section~\ref{sec:(((11)))_2,22}):
\begin{equation}
\frac{d}{dx}Y=\left( \frac{BC}{x}+A_1+A_0 x \right) Y,
\end{equation}
where
\begin{equation*}
B=
\begin{pmatrix}
-Q \\
I
\end{pmatrix}, \quad
C=
\begin{pmatrix}
-P & -PQ+\theta^0
\end{pmatrix}.
\end{equation*}
Note that this is rewritten into the following form:
\begin{equation}
\begin{pmatrix}
\frac{d}{dx}-A_1-A_0 x & -B \\
-C & xI_2
\end{pmatrix}
\begin{pmatrix}
Y \\
Z
\end{pmatrix}=0.
\end{equation}
Applying the Laplace transform, we have
\begin{equation}\label{eq:Laplace(((11)))_2,22}
\left\{
\begin{aligned}
A_0 \hat{Y}'&=-(\xi-A_1)\hat{Y}+B\hat{Z}, \\
\hat{Z}'&=-C\hat{Y}.
\end{aligned}
\right.
\end{equation}
Now we partition the dependent variable $\hat{Y}$ as
\begin{equation*}
\hat{Y}=
\begin{pmatrix}
Y_1 \\
Y_2
\end{pmatrix}
\end{equation*}
where $Y_1$ and $Y_2$ are $2 \times 1$ matrices.
By using (\ref{eq:Laplace(((11)))_2,22}) we can eliminate $Y_1$.
Then the system satisfied by $\begin{pmatrix} Y_2 \\ \hat{Z} \end{pmatrix}$ is
\begin{align*}
\frac{d}{d\xi}
\begin{pmatrix} Y_2 \\ \hat{Z} \end{pmatrix}=
\left[
\xi^2
\begin{pmatrix}
-I & O \\
O & O
\end{pmatrix}
+\xi
\begin{pmatrix}
O & I \\
P & O
\end{pmatrix}
+
\begin{pmatrix}
P-t & -Q \\
PQ-\theta^0 & -P
\end{pmatrix}
\right]
\begin{pmatrix} Y_2 \\ \hat{Z} \end{pmatrix}.
\end{align*}
This system has spectral type $(((2)))(((11)))$.

As to the system $(((2)))_2, 211$ (Section \ref{sec:(((2)))_2,211}), we have to change the independent and dependent variables as
$x \to 1/x,\ Y \to x^{\theta^\infty_1} Y$, so that the Riemann scheme is of the form
\[
\left(
\begin{array}{cc}
  x=0  & x=\infty \ \left( \frac12 \right) \\
\begin{array}{c} 0 \\ 0 \\ \theta^{\infty}_2-\theta^\infty_1 \\ \theta^{\infty}_3-\theta^\infty_1 \end{array} &
\overbrace{\begin{array}{cccc}
   1& 0 & -t/2  & \theta^\infty_1 \\
   1& 0 & -t/2 & \theta^\infty_1 \\
  -1& 0 & t/2 & \theta^\infty_1 \\
  -1& 0 & t/2 & \theta^\infty_1
        \end{array}}
\end{array}
\right).
\]
Then, in the same way as above, we can see the correspondence $(((2)))_2, 211 \leftrightarrow (((2)))(((11)))$.
As the result, we obtain the following correspondences:
\begin{align*}
H^{\mathrm{Mat}}_{\mathrm{II}}:\ 
\stackrel{\infty}{(((11)))}_2, 22 \leftrightarrow (((2)))(((11))), \quad
\stackrel{\infty}{(((2)))}_2, 211 \leftrightarrow (((2)))(((11))).
\end{align*}

\appendix
\section{Data on degenerations}\label{sec:appendix_data}
In this appendix, we give explicit transformations which are used in the calculations of the degenerations.

We note that the gauge parametrizations of the system $22,22,22,211$ and its unramified degenerations,
that is, $(2)(11), 22, 22$, $(2)(2), 22, 211$, $((2))((2)), 211$, $((2))((11)), 22$, $(2)(2), (2)(11)$, and $(((2)))(((11)))$, are different from those in \cite{KNS}.
For example, the gauge parameters of the linear system $22,22,22,211$ in \cite{KNS} is taken to be $U \oplus \mathrm{diag}(v, 1)$.
In this paper, we separate the gauge parameters as
\begin{equation*}
U \oplus \mathrm{diag}(v, 1)=(\mathrm{diag}(v, 1) \oplus \mathrm{diag}(v, 1))\cdot(\mathrm{diag}(v, 1)^{-1} \cdot U \oplus I)
\end{equation*}
and let $v$ be in the matrices $Q$ and $P$ (hence $1/v$ coincides with $u$, see the off-diagonal entries of (\ref{eq:canonical_var})).
We write $\mathrm{diag}(v, 1)^{-1} \cdot U$ as $U$ again.
We have adopted the similar gauge parametrizations as to the other linear systems.

\vspace{3mm}

\noindent
{\bf 2+1+1 $\to$ 3/2+1+1}

\noindent
$(2)(2), 22, 211 \to (2)_2, 22, 211$
\begin{align*}
&\theta^1 = -2\varepsilon^{-1},\ \theta^\infty_1 = \tilde{\theta}^\infty_1+\varepsilon^{-1},\ 
\theta^\infty_2 = \tilde{\theta}^\infty_2+\varepsilon^{-1},\ \theta^\infty_3 = \tilde{\theta}^\infty_3+\varepsilon^{-1},\\
&Q = -(\varepsilon \tilde{t})^{-1}(\tilde{Q}+\theta^0 \tilde{P}^{-1}),\ P = \varepsilon \tilde{t} (1-\tilde{P}),\ 
t = -\varepsilon \tilde{t},\  H = \varepsilon^{-1}\left( -\tilde{H}+\mathrm{tr}\frac{(\tilde{P}-1)(\tilde{Q}+\theta^0\tilde{P}^{-1})}{\tilde{t}} \right),\\
&Y = (x-1)^{-\varepsilon^{-1}}
\begin{pmatrix}
I & O \\
O & -\varepsilon^{-1}
\end{pmatrix}
\tilde{Y}.
\end{align*}

\noindent
$(2)(11),22,22 \to (11)_2,22,22$
\begin{align*}
&\theta^\infty_1 = \varepsilon^{-1},\ \theta^\infty_2 = \tilde{\theta}^\infty_2-\varepsilon^{-1},\ 
\theta^\infty_3 = \tilde{\theta}^\infty_3-\varepsilon^{-1},\\
&Q = \tilde{P},\ P = -\tilde{Q},\ t = \varepsilon \tilde{t},\ H = \varepsilon^{-1}\tilde{H},\\
&Y = \tilde{t}^{-\varepsilon^{-1}}\begin{pmatrix} I & \\ & U \end{pmatrix}^{-1}
\begin{pmatrix} I & \varepsilon^{-1}\tilde{t}^{-1} \\ O & \varepsilon I \end{pmatrix}
\begin{pmatrix} I & -\tilde{P}\tilde{Q}+\theta^0 \\ O & \tilde{t}I \end{pmatrix}\tilde{Y}.
\end{align*}

\noindent
{\bf 3+1 $\to$ 5/2+1}

\noindent
$((2))((2)), 211 \to (((2)))_2, 211$
\begin{align*}
&\theta^0=2\varepsilon^{-6}, \ \theta^\infty_1=\tilde{\theta}^\infty_1-\varepsilon^{-6},\ 
\theta^\infty_2=\tilde{\theta}^\infty_2-\varepsilon^{-6},\ \theta^\infty_3=\tilde{\theta}^\infty_3-\varepsilon^{-6}, \\
&Q=\varepsilon^{-3}+\varepsilon^{-1}\tilde{Q}+\varepsilon(\tilde{P}-\tilde{t}),\ 
P=\varepsilon \tilde{P},\ t=\varepsilon \tilde{t}-2\varepsilon^{-3},\ H=\varepsilon^{-1}\tilde{H}-\varepsilon\,\mathrm{tr}\tilde{P},\\
&x=\varepsilon^{-1}\tilde{x},\ 
Y=\tilde{x}^{\varepsilon^{-6}}\exp\left( \frac{\varepsilon^{-2}}{\tilde{x}} \right)
\begin{pmatrix}
-I & O \\
O & \varepsilon^{-2}I
\end{pmatrix}\tilde{Y}.
\end{align*}

\noindent
$((2))((11)), 22 \to (((11)))_2, 22$
\begin{align*}
&\theta^\infty_1 = -\varepsilon^{-6},\ \theta^\infty_2 = \tilde{\theta}^\infty_2+\varepsilon^{-6},\ 
\theta^\infty_3 = \tilde{\theta}^\infty_3+\varepsilon^{-6},\\
&Q = -\varepsilon \tilde{P},\ P = \varepsilon^{-1}\tilde{Q}-\varepsilon^{-3},\ t = -2\varepsilon^{-3}+\varepsilon \tilde{t},\ H = \varepsilon^{-1}\tilde{H}, \\
&x = \varepsilon \tilde{x},\ Y = e^{\varepsilon^{-2}(\tilde{t}-\tilde{x})}\begin{pmatrix} I & O \\ O & U \end{pmatrix}^{-1}
\begin{pmatrix} I & -\varepsilon^{-3}I \\ O & \varepsilon^{-1}I \end{pmatrix}^{-1}\tilde{Y}.
\end{align*}

\noindent
{\bf 3/2+1+1 $\to$ 5/2+1}

\noindent
$(2)_2, 22, 211 \to (((2)))_2, 211$
\begin{align*}
&\theta^0 = -2\varepsilon^{-3},\ \theta^\infty_1 = \tilde{\theta}^\infty_1+\varepsilon^{-3},\ 
\theta^\infty_2 = \tilde{\theta}^\infty_2+\varepsilon^{-3},\ \theta^\infty_3 = \tilde{\theta}^\infty_3+\varepsilon^{-3},\\
&Q = \varepsilon^{-3}-\varepsilon^{-2}\tilde{Q},\ P = 1-\varepsilon^2\tilde{P},\ t = \varepsilon^{-4}\tilde{t}+\varepsilon^{-6},\ 
H = \varepsilon^4 \tilde{H},\\
&x = -\varepsilon^{-2}\tilde{x},\ Y = \tilde{x}^{-\varepsilon^{-3}}
\begin{pmatrix}
I & O \\
O & -\varepsilon^2 I
\end{pmatrix}\tilde{Y}.
\end{align*}

\noindent
$(2)_2, 22, 211 \to (((11)))_2, 22$
\begin{align*}
&\theta^\infty_1 = \varepsilon^{-3},\ \theta^\infty_2 = \tilde{\theta}^\infty_2-\varepsilon^{-3},
\ \theta^\infty_3 = \tilde{\theta}^\infty_3-\varepsilon^{-3},\\
&Q = -\varepsilon^{-3}+\varepsilon^{-2}(\tilde{Q}-\theta^0 \tilde{P}^{-1}),\ 
P = \varepsilon^2 \tilde{P},\ t = -\varepsilon^{-4}\tilde{t}-\varepsilon^{-6},\ H = -\varepsilon^4 \tilde{H},\\
&x = 1+\frac{1}{\varepsilon^2 \tilde{x}-1},\ 
Y = (1+\varepsilon^2 \tilde{t})^{\varepsilon^{-3}}
\begin{pmatrix}
U & O \\
O & I
\end{pmatrix}^{-1}
G_1
\begin{pmatrix}
I & -\varepsilon^{-1}I \\
O & \varepsilon^2I
\end{pmatrix}
\tilde{Y}.
\end{align*}

\noindent
$(11)_2,22,22 \to (((11)))_2, 22$
\begin{align*}
&\theta^0 = 2\varepsilon^{-3},\ \theta^1 = \tilde{\theta}^0,\ \theta^\infty_2 = \tilde{\theta}^\infty_2-2\varepsilon^{-3},\ 
\theta^\infty_3 = \tilde{\theta}^\infty_3-2\varepsilon^{-3},\\
&Q = -\varepsilon^{-3}+\varepsilon^{-2}(\tilde{Q}-\tilde{\theta}^0 \tilde{P}^{-1}),\ P = \varepsilon^2 \tilde{P},\ 
t = -\varepsilon^{-4}\tilde{t}-\varepsilon^{-6},\ H = -\varepsilon^4 \tilde{H}, \\
&x = 1-\varepsilon^2 \tilde{x},\ Y = (1-\varepsilon^2 \tilde{x})^{-\varepsilon^{-3}}(1+\varepsilon^2 \tilde{t})^{-\varepsilon^{-3}}
\begin{pmatrix} I & \varepsilon^{-1}I \\ O & -\varepsilon^2 I \end{pmatrix}\tilde{Y}.
\end{align*}

\noindent
{\bf 3/2+1+1 $\to$ 3/2+2}

\noindent
$(2)_2, 22, 211 \to (2)_2, (2)(11)$
\begin{align*}
&\theta^0=-\varepsilon^{-1},\ \theta^\infty_2=\tilde{\theta}^\infty_2+\varepsilon^{-1},\ 
\theta^\infty_3=\tilde{\theta}^\infty_3+\varepsilon^{-1},\\
&Q=\varepsilon \tilde{t}\tilde{P},\ P=-(\varepsilon \tilde{t})^{-1}\tilde{Q},\ t=-\varepsilon\tilde{t},\ 
H=-\varepsilon^{-1}\tilde{H}+(\varepsilon\tilde{t})^{-1}\mathrm{tr}(\tilde{P}\tilde{Q}),\\
&x=1-\varepsilon\tilde{x},\ Y=
\begin{pmatrix}
U & O \\
O & I
\end{pmatrix}^{-1}
\begin{pmatrix}
I & -\varepsilon tP \\
-\frac{1}{t}(Z+Q) & \varepsilon(ZP+QP+\theta^0)
\end{pmatrix}
\begin{pmatrix}
\tilde{U} & O \\
O & I
\end{pmatrix}\tilde{Y},
\end{align*}
where $\tilde{U}$ satisfies the gauge equation of $(2)_2,(2)(11)$
\[
\tilde{t}\frac{d\tilde{U}}{d\tilde{t}}\tilde{U}^{-1}=2(\tilde{Q}\tilde{P}+\theta^\infty_1).
\]

\noindent
$(11)_2,22,22 \to (11)_2,(2)(2)$
\begin{align*}
&\theta^0 = -\varepsilon^{-1},\ \theta^1 = \tilde{\theta}^0+\varepsilon^{-1},\\
&Q = \varepsilon \tilde{Q}-(\tilde{\theta}^0 \varepsilon+1)\tilde{P}^{-1},\ P = \varepsilon^{-1}\tilde{P},\ 
t = \varepsilon \tilde{t},\ H = \varepsilon^{-1}\tilde{H}, \\
&x = (\varepsilon \tilde{t})^{-1}\tilde{x},\ Y = \begin{pmatrix} u_1 & 0 \\ 0 & u_2 \end{pmatrix}^{\oplus 2}\tilde{Y},
\end{align*}
where $u_1$ and $u_2$ satisfy
\[
\frac{1}{u_1}\frac{du_1}{dt}=\frac{\theta^\infty_2}{t},\quad
\frac{1}{u_2}\frac{du_2}{dt}=\frac{\theta^\infty_3}{t}.
\]

\noindent
{\bf 2+2 $\to$ 3/2+2}

\noindent
$(2)(2), (2)(11) \to (2)_2, (2)(11)$
\begin{align*}
&\theta^0 = -2\varepsilon^{-1},\ \theta^\infty_1 = \tilde{\theta}^\infty_1+\varepsilon^{-1},\ 
\theta^\infty_2 = \tilde{\theta}^\infty_2+\varepsilon^{-1},\ \theta^\infty_3 = \tilde{\theta}^\infty_3+\varepsilon^{-1},\\
&Q = \varepsilon \tilde{Q},\ P = \varepsilon^{-1}\tilde{P},\ t = \varepsilon \tilde{t},\ H = \varepsilon^{-1} \tilde{H},\\
&Y = x^{-\varepsilon^{-1}}\begin{pmatrix} I & O \\ O & \varepsilon^{-1}I \end{pmatrix}\tilde{Y}.
\end{align*}

\noindent
$(2)(2), (2)(11) \to (2)(2), (11)_2$
\begin{align*}
&\theta^\infty_1=\varepsilon^{-1},\ \theta^\infty_2=\tilde{\theta}^\infty_2-\varepsilon^{-1},\ 
\theta^\infty_3=\tilde{\theta}^\infty_3-\varepsilon^{-1},\\
&Q=\varepsilon \tilde{Q}-(\theta^0\varepsilon+1)\tilde{P}^{-1},\ 
P=\varepsilon^{-1}\tilde{P},\ t=\varepsilon\tilde{t},\ H=\varepsilon^{-1}\tilde{H},\\
&x=\varepsilon \tilde{x},\ Y=\begin{pmatrix} U & O \\ O & I \end{pmatrix}^{-1}
\begin{pmatrix} -\varepsilon P^{-1} & P^{-1} \\ O & I \end{pmatrix}\tilde{Y}.
\end{align*}

\noindent
{\bf 5/2+1 $\to$ 7/2}

\noindent
$(((2)))_2, 211 \to (((((11)))))_2$
\begin{align*}
&\theta^\infty_1=-\varepsilon^{-15}+\frac12,\ \theta^\infty_2=\tilde{\theta}^\infty_2+\varepsilon^{-15}-\frac12,\ 
\theta^\infty_3=\tilde{\theta}^\infty_3+\varepsilon^{-15}-\frac12, \\
&Q=\varepsilon \tilde{Q}+\varepsilon^{-5},\ 
P=\varepsilon^2\frac{\tilde{Q}^2+\tilde{t}}{2}+\varepsilon^{-1}\tilde{P}+\varepsilon^{-4}\tilde{Q}-\varepsilon^{-10},\ 
t=\varepsilon^2 \tilde{t}-3\varepsilon^{-10},\ H=\varepsilon^{-2}\tilde{H}-\frac{\varepsilon}{2}\mathrm{tr}\,\tilde{Q}, \\
&x=-\varepsilon^{16}(\tilde{x}+\varepsilon^{-6}),\ 
Y=\exp (\varepsilon^{-3}\tilde{t})
\begin{pmatrix}
U & O \\
O & I
\end{pmatrix}^{-1}
G_0
\begin{pmatrix}
I & \varepsilon^2 Q \\
O & \varepsilon^2 I
\end{pmatrix}
\tilde{Y}.
\end{align*}

\noindent
$(((11)))_2, 22 \to (((((11)))))_2$
\begin{align*}
&\theta^0 = -2\varepsilon^{-15},\ \theta^\infty_2 = \tilde{\theta}^\infty_2+2\varepsilon^{-15},\ 
\theta^\infty_3 = \tilde{\theta}^\infty_3+2\varepsilon^{-15},\\
&Q = \varepsilon \tilde{Q}+\varepsilon^{-5},\ 
P = \varepsilon^2\frac{\tilde{Q}^2+\tilde{t}}{2}+\varepsilon^{-1}\tilde{P}+\varepsilon^{-4}\tilde{Q}-\varepsilon^{-10},\ 
t = \varepsilon^2 \tilde{t}-3\varepsilon^{-10},\ H \to \varepsilon^{-2}\tilde{H}-\frac{\varepsilon}{2}\mathrm{tr}\,\tilde{Q}, \\
&x = \varepsilon^{-4}\tilde{x}-\varepsilon^{-10},\ Y = (1-\varepsilon^6 \tilde{x})^{-\varepsilon^{-15}}
\begin{pmatrix} I & O \\ O & -\varepsilon^2 I  \end{pmatrix}\tilde{Y}.
\end{align*}

\noindent
{\bf 4 $\to$ 7/2}

\noindent
$(((2)))(((11))) \to (((((11)))))_2$
\begin{align*}
&\theta^\infty_1 = 2\varepsilon^{-15},\ \theta^\infty_2 = \tilde{\theta}^\infty_2-2\varepsilon^{-15},\ 
\theta^\infty_3 = \tilde{\theta}^\infty_3-2\varepsilon^{-15},\\
&Q = \varepsilon \tilde{Q}+\varepsilon^{-5},\ 
P = \varepsilon^2\frac{\tilde{Q}^2+\tilde{t}}{2}+\varepsilon^{-1}\tilde{P}+\varepsilon^{-4}\tilde{Q}-\varepsilon^{-10},\ 
t = \varepsilon^2 \tilde{t}-3\varepsilon^{-10},\ H = \varepsilon^{-2}\tilde{H}-\frac{\varepsilon}{2}\mathrm{tr}\,\tilde{Q}, \\
&x = \varepsilon \tilde{x}+\varepsilon^{-5},\ 
Y = \exp\left( \frac{\varepsilon^{-3}}{2}\tilde{x}^2-\varepsilon^{-9}\tilde{x} \right)
\begin{pmatrix} U & O \\ O & I \end{pmatrix}^{-1}
\begin{pmatrix} I & -\varepsilon^{-3}I \\ O & -\varepsilon^{-8}I \end{pmatrix}\tilde{Y}.
\end{align*}

\noindent
{\bf 3/2+2 $\to$ 7/2}

\noindent
$(2)_2, (2)(11) \to (((((11)))))_2$
\begin{align*}
&\theta^\infty_1=\frac12(1-3\varepsilon^{-5}),\ 
\theta^\infty_2=\tilde{\theta}^\infty_2-\frac12(1-3\varepsilon^{-5}),\ 
\theta^\infty_3=\tilde{\theta}^\infty_3-\frac12(1-3\varepsilon^{-5}), \\
&Q=\varepsilon^{-10}(1-\varepsilon^2 \tilde{Q}),\ 
P=-\varepsilon^8\tilde{P}+\left( \frac{3}{2}\varepsilon^{-5}-1 \right)\varepsilon^{10}(1-\varepsilon^2\tilde{Q})^{-1}
-\varepsilon^5(1+\varepsilon^4\tilde{t})(1-\varepsilon^2\tilde{Q})^{-2}, \\
&t=2\varepsilon^{-15}(1+\varepsilon^4\tilde{t}),\ 
H=\frac{\varepsilon^{11}}{2}\tilde{H}-\mathrm{tr}\left[ \frac{\varepsilon^{10}}{2}(1-\varepsilon^2\tilde{Q})^{-1} \right],\\
&x=2\varepsilon^{-5}(1+\varepsilon^{2}\tilde{x}),\ 
Y=
e^{\varepsilon^{-3}\tilde{x}-\varepsilon^{-1}\tilde{t}}
\begin{pmatrix}
U & O \\
O & I
\end{pmatrix}^{-1}
\begin{pmatrix}
I & -\varepsilon^{-1}I \\
P & -\varepsilon^{-1}P+\frac{2}{t}\varepsilon^{-11}
\end{pmatrix}
\tilde{Y}.
\end{align*}

\noindent
$(2)(2), (11)_2 \to (((((11)))))_2$
\begin{align*}
&\theta^0=3\varepsilon^{-5}-1,\ \theta^\infty_2=\tilde{\theta}^\infty_2-3\varepsilon^{-5}+1,\ 
\theta^\infty_3=\tilde{\theta}^\infty_3-3\varepsilon^{-5}+1,\\
&Q=\varepsilon^{-10}(1-\varepsilon^2 \tilde{Q}),\ 
P=-\varepsilon^8\tilde{P}+\left( \frac{3}{2}\varepsilon^{-5}-1 \right)\varepsilon^{10}(1-\varepsilon^2\tilde{Q})^{-1}
-\varepsilon^5(1+\varepsilon^4\tilde{t})(1-\varepsilon^2\tilde{Q})^{-2}, \\
&t=2\varepsilon^{-15}(1+\varepsilon^4\tilde{t}),\ 
H=\frac{\varepsilon^{11}}{2}\tilde{H}-\mathrm{tr}\left[ \frac{\varepsilon^{10}}{2}(1-\varepsilon^2\tilde{Q})^{-1} \right], \\
&x=\varepsilon^{-10}(\varepsilon^2 \tilde{x}-1),\ 
Y=
\begin{pmatrix}
I & -\varepsilon^{-1}I \\
0 & \varepsilon^4 I
\end{pmatrix}
(\varepsilon^2\tilde{x}-1)^{\frac32 \varepsilon^{-5}}\exp
\left( \varepsilon^{-1}\tilde{t}-\frac{\varepsilon^{-5}}{\varepsilon^2 \tilde{x}-1} \right)\tilde{Y}.
\end{align*}

\noindent
{\bf 3/2+2 $\to$ 3/2+3/2}

\noindent
$(2)_2, (2)(11) \to (2)_2, (11)_2$
\begin{align*}
&\theta^\infty_1=\frac12-\varepsilon^{-1},\ 
\theta^\infty_2=\tilde{\theta}^\infty_2-\frac12+\varepsilon^{-1},\ 
\theta^\infty_3=\tilde{\theta}^\infty_3-\frac12+\varepsilon^{-1}, \\
&Q=-\tilde{Q}(\varepsilon\tilde{P}\tilde{Q}+I),\ 
P=-\varepsilon^{-1}\tilde{Q}^{-1},\ t=\varepsilon\tilde{t},\ H=\varepsilon^{-1}\tilde{H}, \\
&x=-\varepsilon \tilde{x},\ 
Y=\tilde{x}^{\frac12}
\begin{pmatrix}
U & O \\
O & I
\end{pmatrix}^{-1}
\begin{pmatrix}
P^{-1} & \varepsilon^{-1}P^{-1} \\
O & \varepsilon^{-1}I
\end{pmatrix}
\tilde{Y}.
\end{align*}

\noindent
$(2)(2), (11)_2 \to (2)_2, (11)_2$
\begin{align*}
&\theta^0 = -1+2\varepsilon^{-1},\ \theta^\infty_2 = \tilde{\theta}^\infty_2+1-2\varepsilon^{-1},
\theta^\infty_3 = \tilde{\theta}^\infty_3+1-2\varepsilon^{-1},\\
&Q = -\tilde{Q}(\varepsilon \tilde{P}\tilde{Q}+1),\ P = -\varepsilon^{-1}\tilde{Q}^{-1},\ t = \varepsilon \tilde{t},\ 
H = \varepsilon^{-1}\tilde{H},\\
&Y = \tilde{x}^{\varepsilon^{-1}}\tilde{Y}.
\end{align*}


\section{Non-abelian description of matrix Painlev\'e systems}\label{sec:non-abel}
In our study of the four-dimensional Painlev\'e-type equations,
we emphasize their formulation in terms of Hamiltonian systems.

However, concerning the matrix Painlev\'e systems, we have another description of them as ``non-abelian Painlev\'e equations''.
Since it seems to be interesting, we give it in this appendix.

$H^{\mathrm{Mat}}_{\mathrm{VI}}$:
{\small
\begin{equation}
\left\{
\begin{aligned}
t(t-1)\frac{dQ}{dt}&=(Q-t)PQ(Q-1)+Q(Q-1)P(Q-t) \\
&\quad+\delta Q(Q-1)+(-2\alpha-\beta-\gamma-\delta)Q(Q-t)+\gamma(Q-1)(Q-t),\\
t(t-1)\frac{dP}{dt}&=-(Q-1)P(Q-t)P-P(Q-t)PQ-PQ(Q-1)P \\
&\quad-\left[\delta \{P(Q-1)+QP\}+(-2\alpha-\beta-\gamma-\delta)\{P(Q-t)+QP\}+\gamma \{P(Q-t)+(Q-1)P\}\right]  \\
&\quad-\alpha(\alpha+\beta).
\end{aligned}
\right.
\end{equation}
}

$H^{\mathrm{Mat}}_{\mathrm{V}}$:
\begin{equation}
\left\{
\begin{aligned}
t\frac{dQ}{dt}&=Q(Q-1)(P+t)+PQ(Q-1)+\beta Q+\gamma,\\
t\frac{dP}{dt}&=-(Q-1)P(P+t)-P(P+t)Q-\beta P+(\alpha+\gamma)t.
\end{aligned}
\right.
\end{equation}

$H^{\mathrm{Mat}}_{\mathrm{IV}}$:
\begin{equation}
\frac{dQ}{dt}=Q(P-Q-t)+PQ+\beta,\quad
\frac{dP}{dt}=P(-P+Q+t)+QP-\alpha.
\end{equation}

$H^{\mathrm{Mat}}_{\mathrm{III}(D_6)}$:
\begin{equation}\label{eq:non_abel_III(D6)}
t\frac{dQ}{dt}=2QPQ-Q^2+\beta Q+t,\quad
t\frac{dP}{dt}=-2PQP+PQ+QP-\beta P+\alpha.
\end{equation}

$H^{\mathrm{Mat}}_{\mathrm{III}(D_7)}$:
\begin{equation}
t\frac{dQ}{dt}=2QPQ+\alpha Q+t,\quad
t\frac{dP}{dt}=-2PQP-\alpha P-1.
\end{equation}

$H^{\mathrm{Mat}}_{\mathrm{III}(D_8)}$:
\begin{equation}
t\frac{dQ}{dt}=2QPQ+Q, \quad
t\frac{dP}{dt}=-2PQP-P+I-tQ^{-2}.
\end{equation}

$H^{\mathrm{Mat}}_{\mathrm{II}}$:
\begin{align}\label{eq:non_abel_II}
\frac{dQ}{dt}=2P-Q^2-t, \quad \frac{dP}{dt}=PQ+QP+\alpha.
\end{align}

$H^{\mathrm{Mat}}_{\mathrm{I}}$:
\begin{equation}\label{eq:non_abel_I}
\frac{dQ}{dt}=2P, \quad \frac{dP}{dt}=3Q^2+t.
\end{equation}
Eliminating $P$ from the equations (\ref{eq:non_abel_II}) and (\ref{eq:non_abel_I}),
we obtain
\begin{equation*}
\frac{d^2 Q}{dt^2}=6Q^2+2t, \quad
\frac{d^2 Q}{dt^2}=2Q^3+2tQ+2\alpha-1.
\end{equation*}
These are first treated (without constraint on $Q$) in \cite{BS}.


\begin{thebibliography}{}
\bibitem{BS}
S. P. Balandin and V. V. Sokolov,
On the Painlev\'e test for non-Abelian equations,
{\it Phys. Lett. A}\ \textbf{246}, Issues 3--4 (1998), 267--272.


\bibitem{B2}
P. Boalch,
Simply-laced isomonodromy systems,
{\it Publ. Math. Inst. Hautes \'Etudes Sci.}\ \textbf{116}, No. 1 (2012), 1--68.

\bibitem{BV}
D. G. Babbitt and V. S. Varadarajan,
Formal reduction theory of meromorphic differential equations: a group theoretic view,
{\it Pacific J. Math.}~\textbf{109}, No. 1 (1983).


\bibitem{F}
R. Fuchs,
\"Uber lineare homogene Differentialgleichungen zweiter Ordnung mit drei im Endlichen gelegene wesentlich 
singul\"aren Stellen,
{\it Math. Ann.}~\textbf{63} (1907), 301--321.

\bibitem{FS1}
K. Fuji and T. Suzuki,
Drinfeld-Sokolov hierarchies of type $A$ and fourth order Painlev\'{e} systems,
{\it Funkcial. Ekvac.}~{\bf 53} (2010), 143--167. 

\bibitem{Gm}
B. Gambier,
\newblock Sur les \'equations diff\'eretielles du second 
ordre et du premier degr\'e dont l'int\'egrale 
g\'en\'erale est \^^ a points critiques fixes,
\newblock {\it Acta Math.}~{\bf 33} (1910), 1--55.

\bibitem{G}
R. Garnier,
Sur des \'equations diff\'erentielles du troisi\`eme ordre
dont l'int\'egrale g\'en\'erale est uniforme 
et sur une classe d'\'equations nouvelles d'ordre sup\'erieur
dont l'int\'egrale g\'en\'erale a ses points critiques fixes,
\textit{Ann. Sci. \'Ec. Norm. Sup\'er.}~{\bf 29} (1912), 1--126.

\bibitem{HF}
Y. Haraoka and G. Filipuk,
Middle convolution and deformation for Fuchsian systems,
{\it J. Lond. Math. Soc. (2)}\ {\bf 76} (2007), 438--450.

\bibitem{Hrd}
J. Harnad,
\newblock Dual isomonodromic deformations and moment maps to loop algebras,
\newblock {\it Commun. Math. Phys.}~{\bf 166} (1994), 337--365.

\bibitem{Huk}
M. Hukuhara,
\newblock Sur les points singuliers des \'equations diff\'erentielles lin\'eaires, II,
\newblock {\it J. Fac. Sci. Hokkaido Univ.}~{\bf 5} (1937), 123--166.



\bibitem{JMU}
M. Jimbo, T. Miwa, and K. Ueno,
Monodromy preserving deformation of linear ordinary differential equations with rational coefficients. I, 
{\it Phys. 2D} (1981), 306 -- 352.

\bibitem{Katz}
N. M. Katz, Rigid Local Systems,
{\it Princeton Univ. Press}, (1996).

\bibitem{K}
H. Kawakami,
Matrix Painlev\'e systems,
{\it J. Math. Phys.}\ {\bf 56} (2015), doi.org/10.1063/1.4914369.

\bibitem{K2}
H. Kawakami,
Four-dimensional Painlev\'e-type equations associated with ramified linear equations II: Sasano systems,
in preparation.

\bibitem{K3}
H. Kawakami,
Four-dimensional Painlev\'e-type equations associated with ramified linear equations III: Garnier systems and FS systems,
in preparation.

\bibitem{KNS}
H. Kawakami, A. Nakamura, and H. Sakai,
Degeneration scheme of 4-dimensional Painlev\'e-type equations, arXiv:1209.3836.


\bibitem{Ki}
H. Kimura,
The degeneration of the two dimensional Garnier system and the polynomial Hamiltonian structure,
{\it Ann. Mat. Pura Appl. (4)}~{\bf 155} (1989), 25--74.


\bibitem{Lev}
A. H. M. Levelt,
\newblock Jordan decomposition for a class of singular differential operators,
\newblock {\it Ark. Mat.}~{\bf 13 (1)} (1975), 1--27.

\bibitem{OKSO}
Y. Ohyama, H. Kawamuko, H. Sakai, and K. Okamoto,
Studies on the Painlev\'e equations V, third Painlev\'e equations of
special type $P_\mathrm{III}(D_7)$ and $P_\mathrm{III}(D_8)$,
\textit{J. Math. Sci. Univ. Tokyo}~\textbf{13} (2006), 145--204.

\bibitem{OO}
Y. Ohyama and S. Okumura,
A coalescent diagram of the Painlev\'e equations from the viewpoint of isomonodromic deformations,
\textit{J. Phys. A}~\textbf{39} no. 39 (2006), 12129--12151.

\bibitem{O1}
K. Okamoto, 
\newblock Sur les feuilletages associ\'es aux \'equation 
du second ordre \`a points crtiques fixes de P. Painlev\'e,
\newblock {\it Jpn. J. Math.}~{\bf 5} (1979), 1--79.

\bibitem{Ok}
K. Okamoto, 
\newblock Studies on the Painlev\'e equations I.  
\newblock {\it Ann. Mat. Pura Appl.}~{\bf CXLVI} (1987), 337--381;
\newblock II. {\it Jpn. J. Math.}~{\bf 13} (1987), 47--76;
\newblock III. {\it Math. Ann.}~{\bf 275} (1986), 221--255;
\newblock IV. {\it Funkcial. Ekvac.}~{\bf 30} 
(1987), 305--332.

\bibitem{Os2005}
T. Oshima,
A quantization of conjugacy classes of matrices,
{\it Adv. Math.}~\textbf{196} (2005), 124--146.


\bibitem{Os}
T. Oshima,
Fractional calculus of Weyl algebra and Fuchsian differential equations,
\textit{MSJ Memoirs}~\textbf{28} (2012).

\bibitem{P} 
P. Painlev\'e,
Sur les \'equations diff\'erentielles du second ordre \`a
points critiques fixes,
\textit{C. R. Math. Acad. Sci. Paris}~\textbf{127} (1898), 945--948. Oeuvres t. III 35--38.

\bibitem{Sak1}
H. Sakai,
Rational surfaces associated with affine root systems and geometry of the Painlev\'e equations,
{\it Commun. Math. Phys.}~\textbf{220} (2001), 165--229.

\bibitem{Sak2}
H. Sakai,
Isomonodromic deformation and 4-dimensional Painlev\'e type equations,
{\it preprint, University of Tokyo, Mathematical Sciences} (2010).

\bibitem{Ss}
Y. Sasano,
Coupled Painleve VI systems in dimension four with affine Weyl group symmetry of type $D\sp {(1)}\sb 6$. II,
{\it  RIMS K$\hat{o}$ky$\hat{u}$roku Bessatsu}~{\bf B5} 
(2008), 137--152.

\bibitem{Sc}
L. Schlesinger,
\"Uber eine Klasse von Differentialsystemen beliebliger Ordnung mit festen kritischen Punkten,
\textit{J. Reine Angew. Math.}~\textbf{141} (1912), 96--145.

\bibitem{Ts}
T. Tsuda,
UC hierarchy and monodromy preserving deformation,
{\it J. Reine Angew. Math.}~\textbf{690} (2014), 1--34.

\bibitem{Tur}
H. Turrittin,
Convergent solutions of ordinary linear homogeneous differential equations in the neighborhood of an
irregular singular point,
{\it Acta Math.}~{\bf 93} (1955), 27--66.

\bibitem{Wa}
W. Wasow,
Asymptotic Expansions for Ordinary Differential Equations,
{\it John Wiley \& Sons, Inc} (1965).
\end{thebibliography}
\end{document}